\documentclass[reqno]{amsart}

\usepackage{amsmath}
\usepackage{amssymb}
\usepackage{amscd}
\usepackage{graphicx}
\usepackage{epsfig}
\usepackage{psfrag}

\input xy
\xyoption{all}

\theoremstyle{definition}
\newtheorem{thm}{Theorem}
\newtheorem{theorem}[thm]{Theorem}

\newtheorem{lemma}[thm]{Lemma}
\newtheorem{proposition}[thm]{Proposition}
\newtheorem{definition}[thm]{Definition}
\newtheorem{remark}[thm]{Remark}
\newtheorem{example}[thm]{Example}

\numberwithin{thm}{section}
\newcommand{\ft}{{\mathfrak{t}}}
\newcommand{\fb}{{\mathfrak{b}}}
\newcommand{\fh}{{\mathfrak{h}}}
\newcommand{\fg}{{\mathfrak{g}}}

\newcommand{\fp}{{\mathfrak{p}}}
\newcommand{\fv}{{\mathfrak{v}}}

\renewcommand{\SS}{{\mathbb{S}}}
\newcommand{\RR}{{\mathbb{R}}}
\newcommand{\R}{\mathbb{R}}
\newcommand{\ZZ}{{\mathbb{Z}}}
\newcommand{\CC}{{\mathbb{C}}}

\title[GKM Fiber Bundles]{Cohomology of GKM Fiber Bundles}
\author[V. Guillemin]{Victor Guillemin}
\address{Department of Mathematics, MIT, Cambridge, MA 02139}
\email{vwg@math.mit.edu}

\author[S. Sabatini]{Silvia Sabatini}
\address{Department of Mathematics, EPFL, Lausanne, Switzerland}
\email{silvia.sabatini@epfl.ch}

\author[C. Zara]{Catalin Zara}
\address{Department of Mathematics, University of Massachusetts
Boston, MA 02125}
\email{czara@math.umb.edu}

\date{April 15, 2011}

\begin{document}

\maketitle

\begin{abstract}
The equivariant cohomology ring of a GKM manifold is
isomorphic to the cohomology ring of its GKM graph.
In this paper we
explore the implications of this fact for equivariant
fiber bundles for which
the total space and the base space are both GKM and
derive a graph
theoretical version of the Leray-Hirsch theorem.
Then we apply
this result to the equivariant cohomology theory
of flag
varieties.
\end{abstract}

\tableofcontents

\section{Introduction}

Let $T$ be an $n-$dimensional torus, and $M$ a compact, connected $T-$manifold.
The equivariant cohomology ring of $M$,
$H_T^*(M; \RR)$, is an $\SS(\ft^*)-$module, where $\SS(\ft^*) = H_T^*(\text{point})$
is the symmetric algebra on $\ft^*$, the dual of the Lie algebra of $T$. If $H_T^*(M)$
is torsion-free, the restriction map
\begin{equation*}
 i^* \colon H_T^*(M) \to H_T^*(M^T)
\end{equation*}
is injective and hence computing $H_T^*(M)$ reduces to computing the image of $H_T^*(M)$ in $H_T^*(M^T)$.
If $M^T$ is finite, then
$$H_T^*(M^T) = \bigoplus_{p\in M^T} \SS(\ft^*)\; ,$$
with one copy of $\SS(\ft^*)$ for each $p \in M^T$. Determining where
$H_T^*(M)$ sits inside this sum is a challenging problem. However,
one class of spaces $M$ with $H^*_T(M)$ torsion-free for which this
problem has a simple and elegant solution is the one introduced by
Goresky-Kottwitz-MacPherson in their seminal paper \cite{GKM}.
These are now known as \emph{GKM spaces}:
an equivariantly formal space $M$ is a GKM space if $M^T$ is finite
and for every codimension one subtorus $T' \subset T$,
the connected components of $M^{T'}$ are either points or 2-spheres.

To each GKM space $M$ we attach a graph $\Gamma=\Gamma_M$ by
decreeing that the points of $M^T$ are
the vertices of $\Gamma$ and the edges of $\Gamma$ are these
two-spheres. If $S$ is one of the edge
two-spheres, then $S^T$ consists of
exactly two $T-$fixed points, $p$ and $q$.
If $M$ has an invariant almost complex or symplectic structure,
then the isotropy representations on tangent spaces at fixed points
are complex representations and their weights are well-defined.
These data determine a map
$$\alpha \colon E_{\Gamma} \to \ZZ_T^*$$
of oriented edges of $\Gamma$ into the weight lattice of $T$.
This map assigns to the edge 2-sphere
$S$ with North pole $p$ the weight of the isotropy
representation of $T$ on the tangent space to $S$
at $p$. The map $\alpha$ is called the \emph{axial function}
of the graph $\Gamma$. We use it to define
a subring $H_{\alpha}^*(\Gamma_M)$ of $H_T^*(M^T)$ as follows.
Let $c$ be an element of $H_T^*(M^T)$,
\emph{i.e.} a function which assigns to each $p\in M^T$ an
element $c(p)$ of
$H_T^*(\mbox{point})=\mathbb{S}(\mathfrak{t}^*)$.
Then $c$ is in $H_{\alpha}^*(\Gamma_M)$ if and
only if for each edge $e$ of $\Gamma_M$ with vertices $p$
and $q$ as end points,
$c(p)\in \mathbb{S}(\mathfrak{t}^*)$ and
$c(q)\in \mathbb{S}(\mathfrak{t}^*)$ have the same image
in $\mathbb{S}(\mathfrak{t}^*)/\alpha_e\mathbb{S}(\mathfrak{t}^*)$.
(Without the invariant almost complex or symplectic structure, the
isotropy representations are only real representations and the weights
are defined only up to sign; however, that does not change the
construction of $H_{\alpha}^*(\Gamma)$.)
A consequence of a Chang-Skjelbred result (\cite{CS}) is that
$H_{\alpha}^*(\Gamma_M)$ is the image of $i^*$, and
therefore there is an isomorphism of rings
\begin{equation}\label{eq:1.6}
 H_T^*(M) \simeq H_{\alpha}^*(\Gamma_M) \; .
\end{equation}
In a companion paper \cite{GSZ} we prove a fiber bundle
generalization of this result. Let $M$
and $B$ be $T-$manifolds and $\pi \colon M \to B$ be a
$T-$equivariant fiber bundle. If $H_T^*(M)$
is torsion free, then the restriction map
$$i^* \colon H_T^*(M) \to H_T^*(\pi^{-1}(B^T))$$
is injective, and if $B^T$ is finite then $H_T^*(\pi^{-1}(B^T))$ is isomorphic
to
\begin{equation}\label{eq:new1.1}
  \bigoplus_{p \in B^T} H_T^*(F_p)
\end{equation}
with $F_p= \pi^{-1}(p)$. We show in \cite{GSZ} that if $B$ is $GKM$, then
the image of $H_T^*(M)$ in
\eqref{eq:new1.1} can be computed by a generalized
version of \eqref{eq:1.6}. Moreover, if the
fiber bundle is balanced (as defined in \cite{GSZ}),
there is a holonomy action of the groupoid
of paths in $\Gamma$ on the sum \eqref{eq:new1.1} and
the elements which are invariant under this
action form an interesting subring of $H_T^*(M)$.

In this paper we will take the analysis of $H_T^*(M)$
one step further by assuming that $M$ is
also GKM. By interpreting this assumption combinatorially
one is led to a combinatorial notion
which is a central topic of this paper, the notion of a
``fiber bundle of a GKM graph
$(\Gamma_1, \alpha_1)$ over a GKM graph $(\Gamma_2, \alpha_2)$,''
and, associated with this,
the notion of a ``holonomy action'' of the groupoid of paths
in $\Gamma_2$ on the ring
$H_{\alpha_1}(\Gamma_1)$. We will explore below the properties
of such fiber bundles and apply
these results to fiber bundles between generalized flag
varieties; \emph{i.e.} fiber bundles of
the form
\begin{equation}\label{eq:K1K2}
\pi \colon G/P_1 \to G/P_2
\end{equation}
where $G$ is a semi-simple Lie group and $P_1$ and $P_2$
are parabolic subgroups.
In particular we will examine in detail the fiber bundle
\begin{equation}
  \label{eq:example}
  \pi \colon \mathcal{F}l(\CC^n) \to \mathcal{G}r_k(\CC^n)\; ,
\end{equation}
of complete flags in $\CC^n$ over the Grassmannian of
$k-$dimensional subspaces of $\CC^n$
and the analogue of this fibration for the classical groups
of type $B_n$, $C_n$, and $D_n$.
For each of these examples we will compute the subring of
invariant classes in $H_T^*(M)$
(those elements which are fixed by the holonomy action of the
paths in $\Gamma_2$) and show
how the generators of this ring are related to the usual basis
of $H_T^*(M)$, given by equivariant
Schubert classes. These results were inspired by and are related
to results of Sabatini and
Tolman. In \cite{ST} they explore the equivariant cohomology
of fiber bundles where the
total space and the base space are more general symplectic manifolds
with Hamiltonian actions. The theory developed in the
present paper can be regarded as a combinatorial version of the
geometrical theory of symplectic fibrations of
coadjoint orbits, studied in \cite{GLS}.

What follows is a brief table of contents for this paper:
In Section \ref{ssec:motivation} we describe some of the
salient features of the fiber
bundle \eqref{eq:example}. In
Sections~\ref{sec:GKM_graphs}-\ref{ssec:GKM_fiber_bundles} we
briefly review the theory of abstract GKM graphs, following
\cite{GZ1} and \cite{GZ}. We then define abstract versions of fibrations and fiber bundles
between GKM graphs which incorporate
these features, and in
Sections~\ref{ssec:3.1}-\ref{ssec:invariant_classes} we show how
to compute the cohomology ring of such graphs. The main ingredient
in this computation is a holonomy action of the group of based loops in the base on the
cohomology of the fiber graph.

In Section~\ref{sec:flags_as_GKM} we apply this theory to generalized flag manifolds, which
have been extensively studied in the combinatorics literature,
but not from the perspective
of this paper. Let $G$ be a semisimple Lie group, $B$ a Borel
subgroup of $G$ and
$P_1 \subset P_2$ parabolic subgroups containing $B$. Building on results of \cite{GHZ}, in
Section~\ref{ssec:4.1} we describe
the GKM graph associated with the space $P_2/P_1$. In Sections~\ref{ssec:4.2}-\ref{ssec:4.3}
we discuss the fibration of GKM graphs associated with the
fibration of $T-$manifolds
\eqref{eq:K1K2} and compute the group of holonomy automorphisms
associated with this fibration.
In Section~\ref{sec:classical} we specialize to the case where
$G$ is one of the four classical
simple Lie group types, $A_n$, $B_n$, $C_n$, or $D_n$, and,
using iterations of fiber bundles,
give explicit constructions of bases of invariant classes.

In Section~\ref{sec:symmetrization} we construct a second
explicit basis of $H_T^*(G/B)$
consisting of classes that are $W-$invariant. These invariant
classes are obtained from
the equivariant Schubert classes by averaging over the action
of the Weyl group. In Theorem~\ref{th:permuted_classes} we give explicit
combinatorial formulas for the decomposition of twisted Schubert classes,
generalizing earlier results of Tymoczko (\cite[Theorem 4.9]{T}) from
twistings by simple reflections
to actions of general Weyl group elements. We then obtain
formulas for the transition matrix between the basis of invariant
classes consisting of symmetrized Schubert classes and the basis
of invariant classes obtained
through the iterated fiber bundle construction. In addition we
obtain an explicit formula for the decomposition of an invariant
class in the basis of equivariant Schubert classes.

We would like to thank Sue Tolman for her role in inspiring this
work, to Ethan Bolker for helpful
comments on an earlier version, to Allen Knutson and
Alex Postnikov for some very illuminating
remarks concerning the definition of the invariant classes
in the flag manifold case, and to several meticulous referees
whose comments and suggestions improved the presentation of
this paper.

\section{GKM Fiber Bundles}
\label{sec:GKM_bundles}

\subsection{Motivating example}
\label{ssec:motivation}

Let $T^n=(S^1)^n$ be the compact torus of dimension $n$,
with Lie algebra $\ft_n = \RR^n$, and
let $\{ x_1, \ldots , x_n\}$ be the basis of $\ft_n^* \simeq \RR^n$
dual to the canonical basis
of $\RR^n$. Let $\{e_1, \ldots, e_n\}$ be the canonical basis of
$\CC^n$. The torus $T^n$ acts
componentwise on $\CC^n$ by
$$
   (t_1,\ldots,t_n) \cdot (z_1,\ldots, z_n) = (t_1 z_1, \ldots, t_n z_n) \; .
$$
This action induces a $T^n-$action on both $M=\mathcal{F}l(\CC^n)$,
the manifold of complete
flags in $\CC^n$, and $B=\mathcal{G}r_{k}(\CC^n)$, the Grassmannian
manifold of $k-$dimensional
subspaces of $\CC^n$. Let $C = \{ (t,\ldots, t) \; |\;  t\in S^1\}$ be the diagonal circle in
$T^n$ and let $T=T^n/C$. Then $C$ acts trivially on the flag manifold and on Grassmannians, and the induced actions of $T$
on $\mathcal{F}l(\CC^n)$ and on $\mathcal{G}r_{k}(\CC^n)$ are effective. Let
\begin{equation}
\label{eq:first_fibration}
  \pi \colon \mathcal{F}l(\CC^n) \to \mathcal{G}r_{k}(\CC^n) \; ,
\end{equation}
be the map that sends each complete flag $V_{\bullet} = (V_1, \ldots, V_n)$ to its $k-$dimensional
component. Then $(M, B, \pi)$ is a $T-$equivariant fiber bundle.

Since flag manifolds and Grassmannians are GKM spaces, their $T-$equivariant cohomology rings are
determined by fixed point data. These data can be nicely organized using the corresponding GKM graphs,
as follows. For a general GKM space $M$ the fixed point set $M^T$ is finite and is the vertex set of
the GKM graph $\Gamma$. If $T' \subset T$ is a codimension one subtorus of $T$, then the connected
components of the set $M^{T'}$ of $T'-$fixed points are either $T-$fixed points or copies of $\CC P^1$
joining two $T-$fixed points. The edges of the graph $\Gamma$ correspond to these $\CC P^{1}$'s, for
all codimension one subtori $T' \subset T$. An edge $e$ corresponding to a connected component of
$M^{T'}$ is labeled by an element $\alpha_e \in \ft^*$ such that $\ft' = \ker{\alpha_e}$. As explained
in the introduction, the equivariant cohomology ring $H_T^*(M)$ can be computed from the GKM graph $(\Gamma, \alpha)$
associated to $M$, and we will give the details of that construction in Section~\ref{ssec:3.1}.

For the flag manifold $\mathcal{F}l(\CC^n)$, the $T-$fixed point set is indexed by $S_n$, the group
of permutations of $[n] = \{ 1, \ldots, n\}$. A permutation $u=u(1) \ldots u(n)$ of $[n]$ indexes the
fixed flag
$$V_{\bullet}^{u} = (V^{u}_1, \ldots ,V_n^{u})\; ,$$
given by $V_k^{u} = \CC e_{u(1)} \oplus \dotsb \oplus \CC e_{u(k)}$, for all $k = 1,\ldots,n$.

The codimension one subtori $T'$ of $T$ for which the fixed point set is not just the set of $T-$fixed
points are the subtori $T_{ij} = \{ t \in T \; | \; t_i = t_j\} = \exp{(\ker(x_i-x_j))}$. For a fixed
flag $V_{\bullet}^u$, the connected component of $\mathcal{F}l(\CC^n)^{T_{ij}}$ that passes through
$V_{\bullet}^u$ also contains the fixed flag $V_{\bullet}^v$, where $v = (i,j)u$ and $(i,j)$ is the
transposition that swaps $i$ and $j$.

\begin{figure}[h]
  \psfrag{x1-x2}{\small{$x_1\!-\!x_2$}}
  \psfrag{x1-x3}{\small{$x_1\!-\!x_3$}}
  \psfrag{x2-x3}{\small{$x_2\!-\!x_3$}}
  \psfrag{123}{\small{123}}
  \psfrag{213}{\small{213}}
  \psfrag{132}{\small{132}}
  \psfrag{231}{\small{231}}
  \psfrag{312}{\small{312}}
  \psfrag{321}{\small{321}}
  \psfrag{1}{\small{1}}
  \psfrag{2}{\small{2}}
  \psfrag{3}{\small{3}}
  \psfrag{(a)}{(a)}
  \psfrag{(b)}{(b)}
  \includegraphics[height=1.7in]{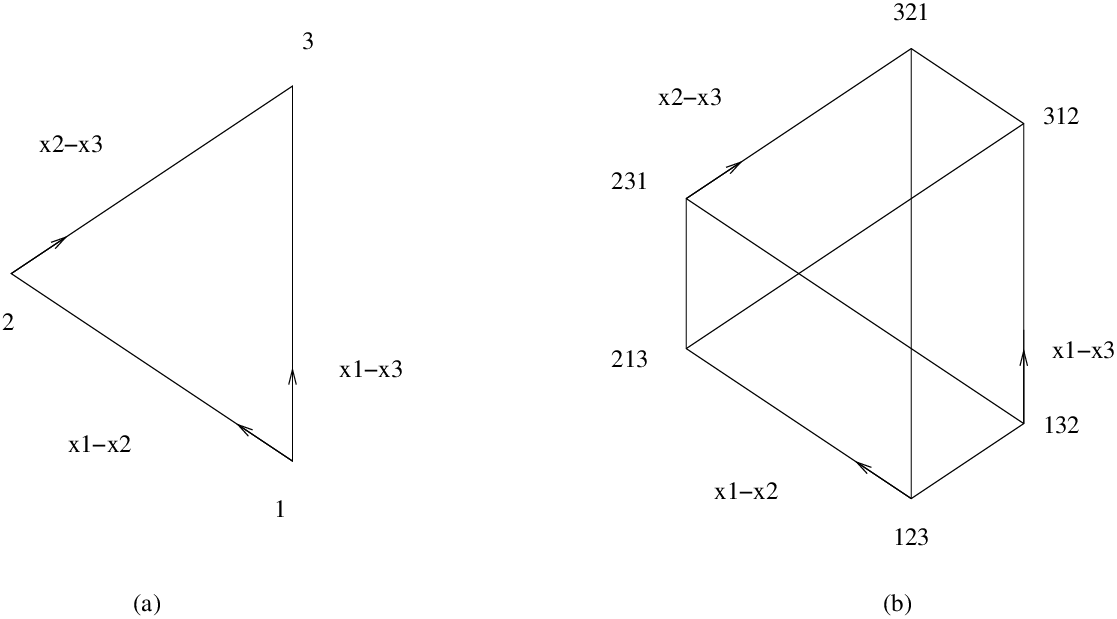}
  \caption{The complete graph $K_3$ (a) and the Cayley graph $(S_3,t)$ (b)}
  \label{fig:k3s3}
\end{figure}

The GKM graph $\Gamma$ of the flag manifold $\mathcal{F}l(\CC^n)$ is the Cayley graph $(S_n,t)$
constructed from the group $S_n$ and generating set $t$, the set of transpositions: the vertices
correspond to permutations in $S_n$ and two vertices are joined by an edge if they differ by a
transposition. If $u \in S_n$, then $u*(i,j) = (u(i),u(j))* u$, so two permutations that differ
by a transposition on the right (operating on \emph{positions}) also differ by a transposition
on the left (operating on \emph{values}). We denote the edge $e$ that joins $u$ and $v=u*(i,j)$ by $u\to v$.
If $1 \leqslant i < j \leqslant n$, then the value of the axial
function $\alpha$ on this edge is
$$\alpha_e = x_{u(i)}-x_{u(j)}\; .$$
We will refer to $\Gamma$ as $S_n$, and
it will be clear from the context when $S_n$ is the graph, the vertex set, or the group of permutations.
Figure~\ref{fig:k3s3}(b) shows the Cayley graph $(S_3,t)$. As a general convention throughout this paper,
edges that are represented by parallel segments have collinear labels. For example,
$\alpha(123,132) = \alpha(231,321) = x_2-x_3$.

For the Grassmannian $\mathcal{G}r_{k}(\CC^n)$, the $T-$fixed point set is indexed by $k-$element
subsets of $[n]$. A subset $I = \{ i_1, \ldots , i_k\}$ corresponds to the fixed $k-$dimensional subspace
$V_I = \CC e_{i_1} \oplus \dotsb \oplus \CC e_{i_k} \; .$
Two vertices are joined by an edge if the intersection of their corresponding $k-$element subsets
is a $(k-1)-$element subset. The resulting graph is the Johnson graph $J(n,k)$. If
$I = (I\cap J ) \cup \{i\}$ and $J = (I\cap J ) \cup \{j\}$, then the value of the axial function
on the edge $e$ from $I$ to $J$ is $\alpha_e = x_i -x_j$. In particular, when $k=1$ we get the
complex projective space $\CC P^{n-1}$, and the associated graph is the  complete graph $K_n$ with
$n$ vertices. The complete graph $K_3$ is shown in Figure~\ref{fig:k3s3}(a).

The discrete version of \eqref{eq:first_fibration} is the morphism of graphs $\pi \colon S_n \to J(n,k)$,
given by $\pi(u) = \{ u(1), \ldots, u(k)\}$. This map is compatible with the axial functions on the two
graphs, and for each vertex $A \in J(n,k)$, the fiber $\pi^{-1}(A)$ is a product $S_k \times S_{n-k}$.
The axial functions on fibers are not identical, but they are compatible in a natural way.

\begin{figure}[h]
  \psfrag{12}{\textbf{12}}
  \psfrag{13}{\textbf{13}}
  \psfrag{14}{\textbf{14}}
  \psfrag{23}{\textbf{23}}
  \psfrag{24}{\textbf{24}}
  \psfrag{34}{\textbf{34}}
  \psfrag{1234}{\scriptsize{1234}}
  \psfrag{1243}{\scriptsize{1243}}
  \psfrag{2134}{\scriptsize{2134}}
  \psfrag{2143}{\scriptsize{2143}}
  \psfrag{1423}{\scriptsize{1423}}
  \psfrag{4123}{\scriptsize{4123}}
  \psfrag{4132}{\scriptsize{4132}}
  \psfrag{4312}{\scriptsize{4312}}
  \psfrag{4321}{\scriptsize{4321}}
  \psfrag{3421}{\scriptsize{3421}}
  \psfrag{3241}{\scriptsize{3241}}
  \psfrag{3214}{\scriptsize{3214}}
  \psfrag{2314}{\scriptsize{2314}}
  \includegraphics[height=2.5in]{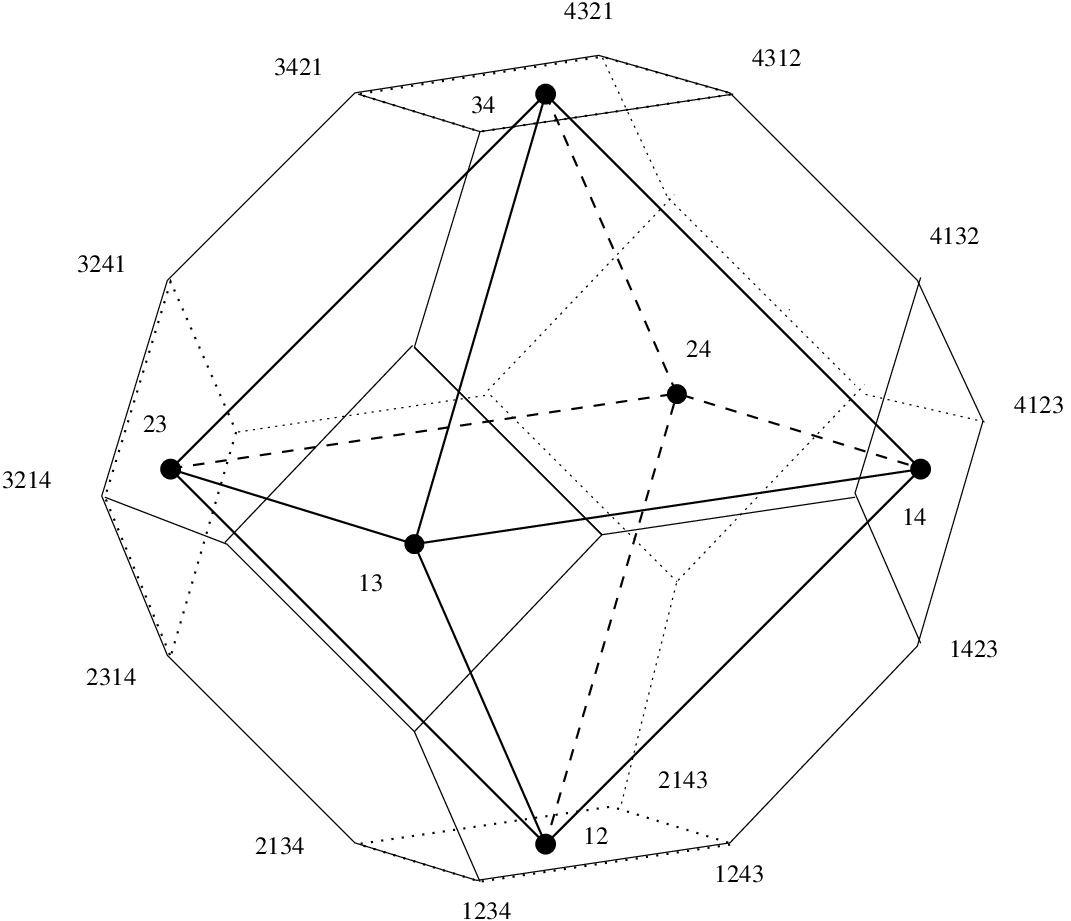}
  \caption{The GKM fiber bundle $S_4 \to J(4,2)$}
  \label{fig:j42}
\end{figure}

The GKM fiber bundle $S_4 \to J(4,2)$  is a combinatorial description of the fiber bundle
$\mathcal{F}l_4(\CC) \to \mathcal{G}r_2(\CC^4)$ that sends a complete flag in $\mathcal{F}l_4(\CC)$
to its two dimensional component. Figure~\ref{fig:j42} shows the graphical representation of this
fiber bundle. The fibers are the squares. (The internal edges of $S_4$ have been omitted.)

This example motivates one of the main goals of this paper: to define the discrete analog of a fiber
bundle between GKM spaces for which the fibers are isomorphic GKM spaces. We then prove a discrete Leray-Hirsch
theorem, showing how one can recover the graph cohomology of the total space from the cohomology of the base and
invariant classes in the cohomology of the fiber.

Then we will revisit the example $\pi \colon \mathcal{F}l(\CC^n) \to \mathcal{G}r_{k}(\CC^n)$ and consider more
general fiber bundles $G/B \to G/P$, with $B \subset P \subset G$ a Borel and parabolic subgroup of a complex
semisimple Lie group $G$, and give a combinatorial description/construction of invariant classes for classical
groups.

\subsection{Abstract GKM Graphs.}
\label{sec:GKM_graphs}
We start by recalling some general definitions (see \cite{GZ1}, \cite{GZ} for more details and motivation).
The reader should have in mind the examples of the Cayley graph $S_n$, the complete graph $K_n$, and the
Johnson graph $J(n,k)$. We will return to these with a summarizing example at the end of
Section~\ref{sec:GKM_bundles}.

Let $\Gamma=(V,E)$ be a regular graph, with $V$ the set of vertices and $E$ the set of oriented edges.
We will consider oriented edges, so each unoriented edge $e$ joining vertices $p$ and $q$ will appear
twice in $E$: once as $(p,q) = p \to q$ and a second time as $(q,p)  = q \to p$. When $e$ is oriented
from $p$ to $q$, we will call $p=i(e)$ the initial vertex of $e$, and $q=t(e)$ the terminal vertex of $e$.
For a vertex $p$, let $E_p$ be the set of oriented edges with initial vertex $p$.

\begin{definition}\label{def:connection}
  Let $e=(p,q)$ be an edge of $\Gamma$, oriented from $p$ to $q$. A \emph{connection along the edge $e$}
  is a bijection $\nabla_e \colon E_p \to E_q$ such that $\nabla_e(p,q) = (q,p)$. A connection on $\Gamma$
  is a family $\nabla = (\nabla_e)_{e\in E}$ of connections along the oriented edges of $\Gamma$, such
  that $\nabla_{(q,p)} = \nabla_{(p,q)}^{-1}$ for every edge $e=(p,q)$ of $\Gamma$.
\end{definition}

\begin{definition}
Let $\nabla$ be a connection on $\Gamma$. A \emph{$\nabla-$compatible axial function} on $\Gamma$ is a
labeling $\alpha \colon E \to \ft^*$ of the oriented edges of $\Gamma$ by elements of a linear
space $\ft^*$, satisfying the following conditions:

\begin{enumerate}
  \item $\alpha(q,p) = - \alpha(p,q)$;
  \item For every vertex $p$, the vectors $\{ \alpha(e) \; | \; e \in E_p\}$ are mutually independent;
  \item For every edge $e=(p,q)$, and for every $e' \in E_p$ we have
  $$ \alpha(\nabla_e(e')) - \alpha(e') = c \alpha(e) \; ,$$
  for some scalar $c \in \RR$ that depends on $e$ and $e'$.
\end{enumerate}
An \emph{axial function} on $\Gamma$ is a labeling $\alpha \colon E \to \ft^*$ that is a
$\nabla-$compatible axial function for some connection $\nabla$ on $\Gamma$.
\end{definition}

\begin{definition}\label{def:GKMgraph}
A \emph{GKM graph} is a pair $(\Gamma, \alpha)$ consisting of a regular graph $\Gamma$ and an
axial function $\alpha \colon E \to \ft^*$ on $\Gamma$.
\end{definition}

\begin{example}[The complete graph]\label{exm:complete}
For the complete graph $\Gamma = K_n$ considered in Section~\ref{ssec:motivation}, the axial function
on \emph{oriented} edges is defined as follows. Let $\ft^*$ be an $n-$dimensional linear space and
$\{x_1,\ldots,x_n\}$ be a basis of $\ft^*$. Define $\alpha \colon E \to \ft^*$ by
$$\alpha(i,j) = x_i-x_j\; .$$
If $\nabla_{(i,j)}  \colon E_i \to E_j$ sends $(i,j)$ to $(j,i)$ and $(i,k)$ to $(j,k)$ for $k \neq i,j$,
then $\nabla$ is a connection compatible with $\alpha$. The image  of $\alpha$ spans the $(n-1)-$dimensional
subspace $\ft_0^*$ generated by $\alpha_1 = x_1-x_2, \ldots, \alpha_{n-1} = x_{n-1} - x_n$.

When $n=2$, the graph $\Gamma$ has two vertices, $1$ and $2$, joined by an edge. The oriented edge from
$1$ to $2$ is labeled $\beta = x_1-x_2$, and the oriented edge from $2$ to $1$ is labeled
$-\beta = x_2-x_1$. The second condition in the definition of an axial function is automatically
satisfied.
\end{example}

\begin{example}[The Cayley graph $(S_n, t)$]\label{exm:permutahedron}
 For the Cayley graph $\Gamma = (S_n,t)$ considered in Section~\ref{ssec:motivation}, the axial
 function on \emph{oriented} edges is defined as follows. Let $\ft^*$ be an $n-$dimensional linear
 space and $\{x_1,\ldots,x_n\}$ be a basis of $\ft^*$. Let $\alpha \colon E \to \ft^*$ be the axial
 function defined as follows. If $u \to v=u(i,j)$ is an oriented edge, with
 $1 \leqslant i < j \leqslant n$, define
$$\alpha(u,v) = x_{u(i)} - x_{u(j)} \; .$$
Note that $\alpha(u,v)$ is determined by the \emph{values} changed from $u$ to $v$.
For an edge $e=u \to v=u(i,j)$, define $\nabla_e \colon E_u \to E_v$ by
\begin{equation}
  \label{eq:connection_permutahedron}
  \nabla_e (u,u(a,b)) = (v,v(a,b))\; .
\end{equation}
Then $\nabla$ is a connection compatible with $\alpha$ and, as above, the image of
$\alpha$ spans the $(n-1)-$dimensional subspace $\ft_0^*$ generated by
$\alpha_1 = x_1-x_2, \ldots, \alpha_{n-1} = x_{n-1} - x_n$.
\end{example}

The examples above show that the image of $\alpha$ may not generate the entire
linear space $\ft^*$. Let $(\Gamma, \alpha)$ be a GKM graph. For a vertex $p$, let
$$\ft_p^* = \text{span}\{\alpha_e \; | \; e\in E_p \} \subset \ft^*$$
be the subspace of $\ft^*$ generated by the image of the axial function on edges
with initial vertex $p$. If $\Gamma$ is connected, then this subspace is the same
for all vertices of $\Gamma$, and we will denote it by $\ft_0^*$. We can co-restrict
the axial function $\alpha \colon E \to \ft^*$ to a function
$\alpha_0 \colon E \to \ft_0^*$, and the resulting pair $(\Gamma, \alpha_0)$ is also a GKM graph.

\begin{definition}
  An axial function $\alpha \colon E \to \ft^*$ is called \emph{effective} if $\ft_0^* = \ft^*$.
\end{definition}

Let $(\Gamma, \alpha)$ be a GKM graph with $\Gamma=(V,E)$ and axial function $\alpha \colon E \to \ft^*$.
Let $\nabla$ be a connection compatible with $\alpha$. Let $\Gamma_0 = (V_0,E_0)$ be a subgraph of $\Gamma$,
with $V_0 \subset V$ and $E_0 \subset E$, such that, if $e \in E$ is an edge with $i(e), t(e) \in V_0$,
then $e \in E_0$.

\begin{definition}
The connected subgraph $\Gamma_0$ is a \emph{$\nabla-$GKM subgraph} if for every edge $e \in E_0$
with $i(e)=p$ and $t(e)=q$, we have $\nabla_e(E_p \cap E_0) = E_q \cap E_0$. The subgraph $\Gamma_0$
is a \emph{GKM subgraph} if it is a  $\nabla-$GKM subgraph for a connection $\nabla$ compatible
with $\alpha$.
\end{definition}

In other words, $\Gamma_0$ is a GKM subgraph if, for some connection $\nabla$ compatible with the axial
function $\alpha$, the connection along edges of $\Gamma_0$ sends edges of $\Gamma_0$ to edges of
$\Gamma_0$ and edges not in $\Gamma_0$ to edges not in $\Gamma_0$. Then the connected subgraph
$\Gamma_0$ is regular, the restriction $\alpha_0$ of $\alpha$ to $E_0$ is an axial function on
$\Gamma_0$, and the connection $\nabla$ induces a connection $\nabla_0$ compatible with $\alpha_0$.
Therefore a GKM subgraph is naturally a GKM graph.

\subsubsection{Isomorphisms of GKM Graphs}
Let $(\Gamma_1,\alpha_1)$ and $(\Gamma_2,\alpha_2)$ be two GKM graphs, with $\Gamma_1=(V_1,E_1)$,
$\alpha_1 \colon E_1 \to \ft_1^*$ and $\Gamma_2=(V_2,E_2)$, $\alpha_2 \colon E_2 \to \ft_2^*$.

\begin{definition}
  An isomorphism of GKM graphs from $(\Gamma_1,\alpha_1)$ to $(\Gamma_2,\alpha_2)$ is a pair
  $(\Phi,\Psi)$, where
  \begin{enumerate}
    \item $\Phi \colon \Gamma_1 \to \Gamma_2$ is an isomorphism of graphs;
    \item $\Psi \colon \ft_1^* \to \ft_2^*$ is an isomorphism of linear spaces;
    \item For every edge $(p,q)$ of $\Gamma_1$ we have
    $$\alpha_2(\Phi(p),\Phi(q)) = \Psi \circ \alpha_1 (p,q) \; . $$
  \end{enumerate}
\end{definition}

The first condition implies that $\Phi$ induces a bijection from $E_1$ to $E_2$, and the third
condition can be restated as saying that the following diagram commutes:
$$
\xymatrix{
E_1 \ar[r]^{\Phi}\ar[d]_{\alpha_1} & E_2 \ar[d]_{\alpha_2} \\
\ft_1^*  \ar[r]^{\Psi}&  \ft_2^* } \,
$$

\subsection{Fiber Bundles of Graphs}

We now introduce special types of morphisms between graphs. Later we will add the GKM package
(axial function and connection) and define the corresponding types of morphisms between GKM graphs.

\subsubsection{Fibrations}
Let $\Gamma$ and $B$ be connected graphs and $\pi \colon \Gamma \to B$ be a morphism of graphs.
By this we mean that $\pi$ is a map from the vertices of $\Gamma$ to the vertices of $B$ such that,
if $(p,q)$ is an edge of $\Gamma$, then either $\pi(p)=\pi(q)$ or else $(\pi(p), \pi(q))$
is an edge of $B$.

When $(p,q)$ is an edge of $\Gamma$ and $\pi(p)=\pi(q)$, we will say that the edge $(p,q)$ is
\emph{vertical}; otherwise $(\pi(p),\pi(q))$ is an edge of $B$ and we will say that $(p,q)$ is
\emph{horizontal}. For a vertex $q$ of $\Gamma$, let $E_q^{\perp}$ be the set of vertical edges
with initial vertex $q$, and let $H_q$ be the set of horizontal edges with initial vertex $q$.
Then $E_q = E_q^{\perp} \cup H_q$ and $\pi$ canonically induces a map
$(d\pi)_q \colon H_q \to (E_B)_{\pi(q)}$ given by
\begin{equation}
  \label{eq:dpi}
   (d\pi)_q(q,q') = (\pi(q),\pi(q'))\; .
\end{equation}

\begin{definition}\label{def:graphfibration}
  The morphism of graphs $\pi \colon \Gamma \to B$ is a
  \emph{fibration of graphs}\footnote{This is what we
called \emph{submersion} in \cite{GZ}. This definition of a fibration of graphs is different
from the one introduced in \cite{BV}. We work with undirected graphs,
and our morphisms of graphs allow
edges to collapse.} if for every vertex $q$ of $\Gamma$, the map
$(d\pi)_q \colon H_q \to (E_B)_{\pi(q)}$ is bijective.
\end{definition}

Fibrations have the \emph{unique lifting of paths} property: Let $\pi \colon \Gamma \to B$
be a fibration, $(p_0,p_1)$ an edge of $B$, and $q_0 \in \pi^{-1}(p_0)$ a point in the fiber
over $p_0$.
Since $(d\pi)_{q_0} \colon H_{q_0} \to (E_B)_{p_0}$ is a bijection, there exists a unique
edge $(q_0,q_1)$ such that $(d\pi)_{q_0} (q_0,q_1) = (p_0,p_1)$. We will say that $(q_0,q_1)$
is \emph{the lift} of $(p_0,p_1)$ at $q_0$. If $\gamma$ is a path $p_0 \to p_1 \to \dotsb \to p_m$
in $B$ and $q_0 \in \pi^{-1}(p_0)$ is a point in the fiber over $p_0$, then we can lift $\gamma$
uniquely to a path $\widetilde{\gamma}(q_0) = q_0 \to q_1 \to \dotsb \to q_m$ in $\Gamma$
starting at $q_0$, by successively lifting the edges of $\gamma$.

\subsubsection{Fiber Bundles}
Let $\pi \colon \Gamma \to B$ be a fibration of graphs.
For a vertex $p$ of $B$, let $V_p=\pi^{-1}(p) \subset V$ and let $\Gamma_{\!p}$ be the induced
subgraph of $\Gamma$ with vertex set $V_p$. For every edge $(p,q)$ of $B$, define a map
$\Phi_{p,q} \colon V_{p} \to V_{q}$ as follows. For $p' \in V_p$, define $\Phi_{p,q}(p') = q'$,
where $(p',q')$ is the lift of $(p,q)$ at $p'$. It is easy to see that $\Phi_{p,q}$ is bijective,
with inverse $\Phi_{q,p}$.  What is not true, in general, is that $\Phi_{p,q}$ is an isomorphism
of graphs from $\Gamma_{\!p}$ to $\Gamma_{\!q}$.

\begin{example}\label{exm:fibration}
  Let $\Gamma$ be the regular $3-$valent graph consisting of two quadrilaterals
  $(p_1,p_2,p_3,p_4)$ and $(q_1,q_3,q_2,q_4)$ joined by edges $(p_i,q_i)$ for i=1,2,3,4.
  (See Figure~\ref{fig:twisted}.)

\begin{figure}[h]
  \psfrag{a1}{$p_1$}
  \psfrag{a2}{$p_2$}
  \psfrag{a3}{$p_3$}
  \psfrag{a4}{$p_4$}
  \psfrag{b1}{$q_1$}
  \psfrag{b2}{$q_2$}
  \psfrag{b3}{$q_3$}
  \psfrag{b4}{$q_4$}
  \psfrag{a}{p}
  \psfrag{b}{q}
  \includegraphics[height=1.5in]{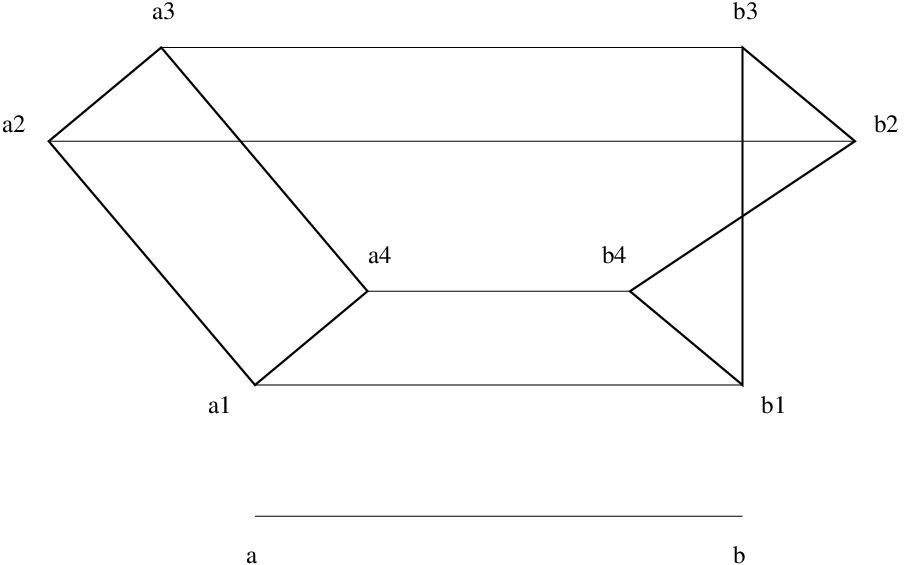}
  \caption{Twisted fibration}
  \label{fig:twisted}
\end{figure}

  Let $B$ be a graph with two vertices $p$ and $q$ joined by an edge. Let $\pi \colon \Gamma \to B$
  be the morphism of graphs $\pi(p_i) = p$ and $\pi(q_i) = q$ for $i=1,2,3,4$. Then $\pi$ is a
  fibration and $\Phi_{p,q}(p_i) =q_i$ for $i=1,2,3,4$. However, $(p_1,p_2)$ is an edge in
  $\Gamma_{\!p}$, but $(q_1,q_2)$ is not an edge in $\Gamma_{\!q}$. While the fibers
  $\Gamma_{\!p}$ and $\Gamma_{\!q}$ are isomorphic as graphs, the map $\Phi_{p,q}$ is not an
  isomorphism.
\end{example}

We will be interested in fibrations for which $\Phi_{p,q}$ is an isomorphism of graphs from
the fiber $\Gamma_{\!p}$ to the fiber $\Gamma_{\!q}$.

\begin{definition}\label{def:fiberbundle}
A fibration $\pi \colon \Gamma \to B$ is a \emph{fiber bundle}\footnote{This is what we
called \emph{fibration} in \cite{GZ}} if for every edge $(p,q)$ of $B$, the map
$\Phi_{p,q} \colon \Gamma_{\!p} \to \Gamma_{\!q}$ is a morphism of graphs.
\end{definition}

If $\pi \colon \Gamma \to B$ is a fiber bundle then $\Phi_{p,q}$ is bijective, and
both $\Phi_{p,q} \colon \Gamma_{\!p} \to \Gamma_{\!q}$ and
$\Phi_{p,q}^{-1} = \Phi_{q,p} \colon \Gamma_{\!q} \to \Gamma_{\!p}$ are morphisms of graphs.
Therefore the maps $\Phi_{p,q}$ are isomorphisms of graphs. The simplest example of a fiber
bundle is the projection of a direct product of graphs onto one of its factors,
$\pi \colon \Gamma =B \times F \to B$. We will call such fiber bundles \emph{trivial bundles}.

\subsection{GKM Fiber Bundles}
\label{ssec:GKM_fiber_bundles}

We now add the GKM package to a fibration, and define GKM fibrations. Let $(\Gamma,\alpha)$ and
$(B,\alpha_B)$ be two GKM graphs, with axial functions $\alpha \colon E \to \ft^*$ and
$\alpha_B \colon E_B \to \ft^*$ taking values in the same linear space $\ft^*$. Let $\nabla$
and $\nabla_B$ be connections on $\Gamma$ and $B$, compatible with $\alpha$ and $\alpha_B$,
respectively.

\begin{definition}\label{def:GKMfibration}
  A map $\pi \colon (\Gamma,\alpha) \to (B,\alpha_B)$ is a \emph{$(\nabla,\nabla_B)-$GKM fibration}
  if it satisfies the following conditions:
\begin{enumerate}
  \item $\pi$ is a fibration of graphs;
  \item If $e$ is an edge of $B$ and $\widetilde{e}$ is any lift of $e$, then
  $\alpha(\widetilde{e}) = \alpha_B(e)$;
  \item Along every edge $e$ of $\Gamma$ the connection $\nabla$ sends horizontal edges
  into horizontal edges and vertical edges into vertical edges;
  \item The restriction of $\nabla$ to horizontal edges is compatible with $\nabla_B$, in
  the following sense:
Let $e=(p,q)$ be an edge of $B$ and $\widetilde{e} = (p',q')$ the lift of $e$ at $p\,'$.
Let $e' \in E_p$ and $e''= (\nabla_B)_{e}(e') \in E_{q}$. If $\widetilde{e'}$ is the lift
of $e'$ at $p\,'$ and $\widetilde{e''}$ is the lift of $e''$ at $q'$ then
$$(\nabla)_{\widetilde{e}} (\widetilde{e'}) = \widetilde{e''} \; .$$
\end{enumerate}
A map $\pi \colon (\Gamma,\alpha) \to (B,\alpha_B)$ is a \emph{GKM fibration} if it is a
$(\nabla,\nabla_B)-$GKM fibration for some connections $\nabla$ and $\nabla_B$ compatible
with $\alpha$ and $\alpha_B$.
\end{definition}

If $\pi \colon (\Gamma,\alpha) \to (B,\alpha_B)$ is a GKM fibration, then for each $p \in B$,
the fiber $(\Gamma_{\!p}, \alpha)$ is a GKM subgraph of $(\Gamma, \alpha)$. Let $\fv_p^*$ be
the subspace of $\ft^*$ generated by values of axial functions $\alpha_e$,
for edges $e$ of $\Gamma_{\! p}$. Then the axial function on $\Gamma_{\!p}$ can be co-restricted
to $\alpha_p$, from the oriented edges of $\Gamma_{\!p}$ to $\fv_p^*$, and $(\Gamma_{\!p}, \alpha_p)$
is a GKM graph.

Suppose now that $\pi$ is both a GKM fibration and a fiber bundle of graphs. Let $e=(p,q)$ be an
edge of $B$. We say that the transition isomorphism $\Phi_{p,q} \colon \Gamma_{\,p} \to \Gamma_{\,q}$
is \emph{compatible with the connection} on $\Gamma$ if for every lift $\tilde{e}=(p_1,q_1)$ of $e$
and for every edge $e' = (p_1,p_2)$ of $\Gamma_{\,p}$, the connection along $\tilde{e}$ moves $e'$
into the edge $e'' = (q_1,q_2) = (\Phi_{p,q}(p_1), \Phi_{p,q}(p_2))$ of $\Gamma_{\,q}$.

\begin{definition}\label{def:GKMfiberbundles}
  A GKM fibration  $\pi \colon (\Gamma,\alpha) \to (B,\alpha_B)$ is a \emph{GKM fiber bundle} if
  $\pi$ is a fiber bundle and for every edge $e=(p,q)$ of $B$:
  \begin{enumerate}
    \item The transition isomorphism $\Phi_{p,q}$ is compatible with the connection of $\Gamma$.
    \item There exists a linear isomorphism $\Psi_{p,q} \colon \fv_p^* \to \fv_q^*$ such that
$$\Upsilon_{p,q} = (\Phi_{p,q}, \Psi_{p,q}) \colon (\Gamma_{\!p},\alpha_p) \to (\Gamma_{\!q}, \alpha_q)$$
is an isomorphism of GKM graphs.
  \end{enumerate}
\end{definition}

For a GKM fiber bundle $\pi \colon (\Gamma,\alpha) \to (B,\alpha_B)$ we can be more specific about
the transition isomorphisms $\Psi_{p,q}$. Let $(p,q)$ be an edge of $B$, let $(p',p'')$ be an edge
of $\Gamma_{\!p}$, and let $(q', q'')$ be the corresponding edge of $\Gamma_{\!q}$. The compatibility
condition along the edge $(p',q')$ implies that $\alpha_{q',q''} - \alpha_{p',p''}$ is a multiple
of $\alpha_{p',q'} = \alpha_{p,q}$, hence
there exists a unique constant $c=c(\alpha_{p',p''})$ such that
$$\Psi_{p,q}(\alpha_{p\,',p\,''} )= \alpha_{p\,',p\,''} + c(\alpha_{p',p''}) \alpha_{p,q} \; .$$
The linearity of $\Psi_{p,q}$ implies that there exists a unique linear function $c \colon \fv_p^* \to \RR$
such that
$$\Psi_{p,q}(x) = x+ c(x)\alpha_{p,q}$$
for all $x \in \fv_p^*$.

For a path $\gamma \colon \; p_0 \to p_1 \to \dotsb \to p_{m-1} \to p_m$ in $B$ from $p_0$ to $p_m$,
let
$$\Upsilon_{\gamma} = \Upsilon_{p_{m-1},p_m} \circ \dotsb \circ \Upsilon_{p_0,p_1} \colon
(\Gamma_{\!p_0},\alpha_{p_0}) \to (\Gamma_{p_m}, \alpha_{p_m})$$
be the GKM graph isomorphism given by the composition of the transition maps. Let $p \in B$ be a
vertex, and let $\Omega(p)$ be the set of all loops in $B$ that start and end at $p$. If
$\gamma \in \Omega(p)$ is a loop based at $p$, then $\Upsilon_{\gamma}$ is an automorphism of
the GKM graph $(\Gamma_{\!p},\alpha_p)$. The \emph{holonomy group} of the fiber $\Gamma_p$ is
the group
$$\text{Hol}_{\pi}(\Gamma_{\!p}) = \{ \Upsilon_{\gamma} \; | \; \gamma \in \Omega(p) \} \leqslant
 \text{Aut}(\Gamma_{\!p}, \alpha_p) \; .$$

If the base $B$ is connected, then all the fibers are isomorphic as GKM graphs. Let $(F,\alpha_F)$
be a GKM graph isomorphic to all fibers, with $\alpha_F \colon E_F \to \ft_F^*$, and, for each
vertex $p$ of $B$, let $\rho_p = (\varphi_p,\psi_p) \colon (F,\alpha_F) \to (\Gamma_{\!p},\alpha)$
be a fixed isomorphism of GKM graphs. For every edge $(p,p\,')$ of $B$, let $\rho_{p,p\,'} =
(\varphi_{p,p\,'},\psi_{p,p\,'}) \colon (F,\alpha_F) \to (F,\alpha_F)$ be the automorphism of
$(F,\alpha_F)$ given by
\begin{align*}
  \varphi_{p,p\,'} & = \varphi_{p\,'}^{-1} \circ \Phi_{p,p\,'} \circ \varphi_p \\
  \psi_{p,p\,'} & = \psi_{p\,'}^{-1} \circ \Psi_{p,p\,'} \circ \psi_p \; .
\end{align*}
If $\gamma$ is any path in $B$, then the composition of the transition maps along the edges of
$\gamma$ defines an automorphism $\rho_{\gamma} = (\varphi_{\gamma}, \psi_{\gamma})$ of
$(F, \alpha_F)$. Let $p$ be a vertex of $B$ and
$$\text{Hol}(F,p) = \{ \rho_{\gamma} \; | \; \gamma \in \Omega(p) \} \subset \text{Aut}(F, \alpha_F) \; .$$
Then $\text{Hol}(F,p)$ is a subgroup of $\text{Aut}(F, \alpha_F)$ and if $p, p\,'$ are
vertices of $B$, then $\text{Hol}(F,p)$ and $\text{Hol}(F,{p\,'})$ are conjugated by
$\rho_{\gamma}$, where $\gamma$ is any path in $B$ connecting $p$ and $p\,'$.

\subsection{Example}
\label{ssec:sum_example}
In this section we return to $\pi \colon \mathcal{F}l(\CC^n) \to \CC P^{n-1}$
(as a particular case of $\mathcal{F}l(\CC^n) \to \mathcal{G}r_k(\CC^n)$).  We show
that the discrete version, $\pi \colon S_n \to K_n$, given by $\pi(u) = u(1)$, is an
abstract GKM fiber bundle.

\subsubsection{$\pi$ is a GKM fibration} Clearly $\pi$ is a morphism of graphs, because
in $K_n$ all vertices are joined by edges.  Moreover, let $u$ and $v=u(i,j)$ (with
$1 \leqslant i<j \leqslant n$) be adjacent vertices in $S_n$. If $i \neq 1$, then $\pi(u) = \pi(v)$,
hence the edge $u \to v$ is vertical. If $i=1$, then $\pi(v) = u(j) \neq u(1)$, hence the edge
$u \to v$ is horizontal.

Let $(d\pi)_u \colon H_u \to E_{\pi(u)}$ be the induced map \eqref{eq:dpi}. If $\tilde{e}$ is
the horizontal edge $u \to v = u(1,j)$, then
$(d\pi)_u (\tilde{e}) = e$, the edge of $K_n$ joining $u(1)$ and $u(j)$. Therefore $(d\pi)_u$
is bijective, hence $\pi$ is a fibration of graphs.

\begin{figure}[h]
  \psfrag{x1-x2}{$x_1\!-\!x_2$}
  \psfrag{x1-x3}{$x_1\!-\!x_3$}
  \psfrag{x2-x3}{$x_2\!-\!x_3$}
  \psfrag{123}{\textbf{1}23}
  \psfrag{213}{\textbf{2}13}
  \psfrag{132}{\textbf{1}32}
  \psfrag{231}{\textbf{2}31}
  \psfrag{312}{\textbf{3}12}
  \psfrag{321}{\textbf{3}21}
  \psfrag{1}{1}
  \psfrag{2}{2}
  \psfrag{3}{3}
  \psfrag{pi}{$\pi$}
  \psfrag{fib1}{\small{$\pi^{-1}(1)$}}
  \psfrag{fib2}{\small{$\pi^{-1}(2)$}}
  \psfrag{fib3}{\small{$\pi^{-1}(3)$}}
  \includegraphics[height=2in]{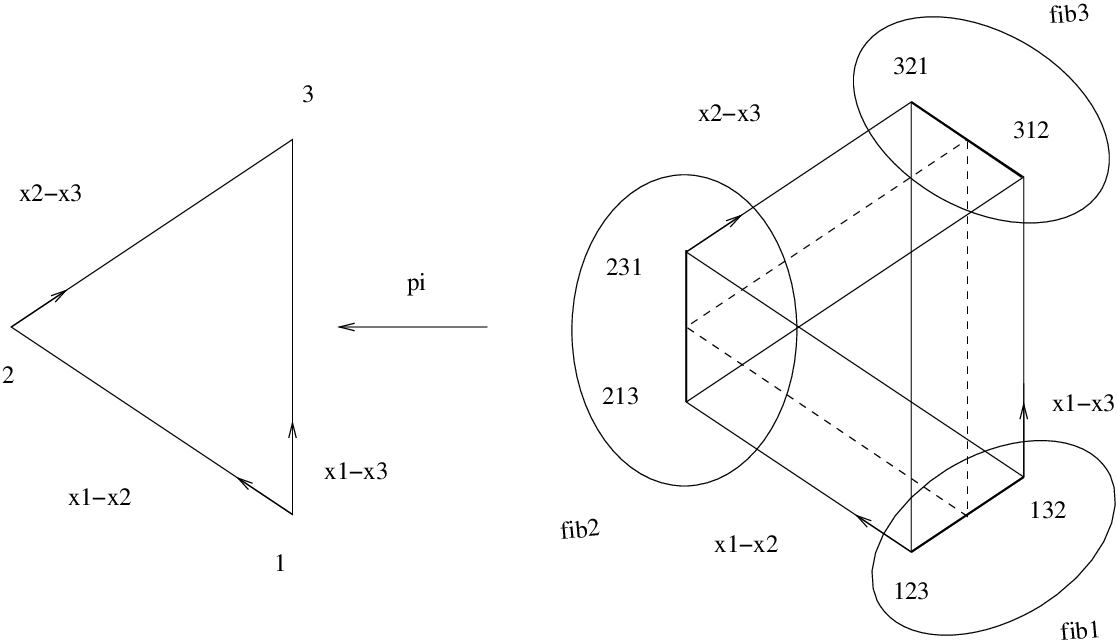}
  \caption{Fibration $S_3 \to K_3$}
  \label{fig:s3k3fibration}
\end{figure}

The case $n=3$ is shown in Figure~\ref{fig:s3k3fibration}. If $\gamma$ is the cycle
$1 \to 2 \to 3 \to 1$ in $K_3$, then the lift of $\gamma$ at $123$ is the path
$\widetilde{\gamma} : 123 \to 213 \to 312 \to 132$ in $S_3$.

In Section~\ref{ssec:4.2} we will prove in a more general case that $\pi$ is a GKM
fibration. For now, we notice that it is compatible with the axial functions: if
$\tilde{e}$ is the horizontal edge in $S_n$ from $u$ to $v = u(1,j)$, then $\tilde{e}$
is a lift of the edge $e$ in $K_n$ from $u(1)$ to $u(j)$, and both edges have the same
label, $x_{u(1)} - x_{u(j)}$.

\subsubsection{Transition isomorphisms} For each $i \in [n]$, the fiber $\Gamma_i = \pi^{-1}(i)$
consists of all permutations $u \in S_n$ for which $u(1) = i$ and is isomorphic, as a graph,
with the Cayley graph $(S_{n-1},t)$. For $1 \leqslant i \neq j \leqslant n$, the transition
map $\Phi_{i,j} \colon \Gamma_i \to \Gamma_j$ is given by $\Phi_{i,j}(u) = (i,j)u$. Let $u$
be a vertex of $\Gamma_i$ and $u\to u' = u(a,b)$ an edge of $\Gamma_{\!i}$, hence
$2 \leqslant a < b \leqslant n$. If $v = \Phi_{i,j}(u)$ and $v'= \Phi_{i,j}(u')$,
then $v' = (i,j)u' = (i,j)u(a,b) = v(a,b)$, hence $v$ and $v'$ are joined by an edge
in $\Gamma_{\!j}$. This shows that $\Phi_{i,j}$'s are morphisms of graphs, and therefore
$\pi$ is a fiber bundle.

The subspace generated by the values of the axial function on the edges of $\Gamma_i$ is
the $(n-1)-$dimensional space
$$\fv_i^* = \text{span}_{\R} \{  x_r -x_s \, | \, 1 \leqslant r \neq i \neq s \leqslant n\}$$
and similarly for $\Gamma_{\!j}$. Let $\alpha_i$ and $\alpha_j$ be the axial functions on
$\Gamma_{\!i}$ and $\Gamma_{\!j}$. Then $\alpha_i(u,u') = x_{u(a)} - x_{u(b)}$, and
$\alpha_j(v,v') = x_{(i,j)u(a)}-x_{(i,j)u(b)}$. Let $\Psi_{i,j}$ be the linear automorphism
of $\ft^*$ induced by $\Psi_{i,j}(x_r) = x_{(i,j)r}$, for $1 \leqslant r \leqslant n$.
Then $\Psi_{i,j}$ induces an isomorphism $\Psi_{i,j} \colon \fv_i^* \to \fv_j^*$ and
$$\alpha_j(\Phi_{i,j}(u),\Phi_{i,j}(u')) = \Psi_{i,j}(\alpha_i(u,u'))\; ,$$
which proves that $(\Phi_{i,j}, \Psi_{i,j}) \colon (\Gamma_{\!i}, \alpha_i) \to
(\Gamma_{\!j}, \alpha_j)$ is an isomorphism of GKM graphs. Since the fibers are canonically
isomorphic as GKM graphs, the map $\pi \colon S_n \to K_n$ is a GKM fiber bundle.

\subsubsection{Typical fiber} For $1 \leqslant i \leqslant n$, the fiber
$(\Gamma_{\!i}, \alpha_i)$ is isomorphic to $S_{n-1}$, and we construct an explicit
isomorphism $\varphi_i \colon S_{n-1} \to \Gamma_{\!i}$. For a permutation $u \in S_{n-1}$,
let $\tilde{u} = u(1)u(2)\dotsb u(n\!-\!1)n \in S_n$. For $1 \leqslant a < b \leqslant n$,
let $c_{a,b}$ be the cycle $a \to a+1 \to \dotsb \to b \to a$, and let $c_{b,a} = c_{a,b}^{-1}$.
Then the map $\varphi_i \colon S_{n-1} \to \Gamma_{\!i}$,
$$\varphi_i(u) = c_{i,n}\tilde{u}c_{n,1}$$
is a graph isomorphism between $S_{n-1}$ and $\Gamma_{\!i}$. The cycle $c_{i,n}$, operating
on values, moves the value $i$ to the last position and preserves the relative order of
the values on the other positions. The cycle $c_{n,1}$, operating on positions, moves $i$
from the last position to the first and then shifts all the other positions to the right by one.

Let $\psi_i$ be the linear isomorphism induced by $\psi_i(x_k) = x_{c_{i,n}(k)}$ for all
$1 \leqslant k \leqslant n$. If $u \in S_{n-1}$ and $v=u(a,b)$, with $1 \leqslant a < b \leqslant n\!-\!1$, then
$$\alpha_i (\varphi_i (u), \varphi_i(v)) = \psi_i(\alpha(u,v)) \; ,$$
hence $(\varphi_i, \psi_i) \colon S_{n-1} \to \Gamma_{\!i}$ is an isomorphism of GKM graphs.

\subsubsection{Holonomy action on the fiber} Let $Hol(\Gamma_n)$ be the holonomy group of
the fiber $\Gamma_n$. It is generated by compositions of transition isomorphisms along loops
in $K_n$ based at $n$. Each such nontrivial loop can be decomposed into triangles
$\gamma_{ij} \colon n \to i \to j \to n$, and for such a triangle we have $(j,n)(i,j)(n,i) = (i,j)$,
hence the corresponding element of $Hol(\Gamma_n)$ generated by $\gamma_{ij}$ is
$$\Upsilon_{\gamma_{ij}} = (\Phi_{\gamma_{ij}}, \Psi_{\gamma_{ij}})\; ,$$
with $\Phi_{\gamma_{ij}} (u) = (i,j)u$ and $\Psi_{\gamma_{ij}}(x_r) = x_{(i,j)r}$.

Since every permutation in $S_{n-1}$ can be decomposed into transpositions, it follows that
$$\text{Hol}(\Gamma_n) = \{ \Upsilon_{w} = (\Phi_w, \Psi_w) \; | \; w \in S_{n-1} \} \simeq S_{n-1}\; , $$
where, for a permutation $w \in S_{n-1}$, $\Phi_w \colon \Gamma_n \to \Gamma_n$ is given by
$\Phi_w(u) = wu$, and $\Psi_w (x_a) = x_{w(a)}$.

Since the holonomy actions are conjugated, it follows that the holonomy group of all fibers
are isomorphic to $S_{n-1}$.

\section{Cohomology of GKM Fiber Bundles}
\label{sec:Cohomology}

Let $\pi \colon (\Gamma, \alpha) \to (B, \alpha_B)$ be a GKM fiber bundle, with typical fiber
$(F, \alpha_F)$. One of the main goals of this paper is to describe how the cohomology ring of
the total space ($\Gamma, \alpha)$ is determined by the cohomology rings of the base
$(B, \alpha_B)$ and the fiber $(F, \alpha_F)$ and the holonomy action of the base on the fiber.
We start by recalling the construction of the cohomology ring of a GKM graph.

\subsection{Cohomology of GKM graphs}
\label{ssec:3.1}

Let $(\Gamma, \alpha)$ be a GKM graph, with $\Gamma = (V,E)$ a regular graph and $\alpha \colon E \to \ft^*$
an axial function. Let $\SS(\ft^*)$ be the symmetric algebra of $\ft^*$; if $\{x_1,\ldots,x_n\}$ is
a basis of $\ft^*$, then $\SS(\ft^*) \simeq \RR[x_1,\ldots,x_n]$.

\begin{definition}\label{def:cohclass}
  A cohomology class on $(\Gamma, \alpha)$ is a map $\omega \colon V \to \SS(\ft^*)$ such that for every
  edge $e=(p,q)$ of $\Gamma$, we have
\begin{equation}\label{eq:compcond}
  \omega(q) \equiv \omega(p) \pmod{\alpha_e} \; .
\end{equation}
\end{definition}

The compatibility condition \eqref{eq:compcond} means that $\omega(q) - \omega(p) = \alpha_e f$,
for some element $f \in \SS(\ft^*)$, and is equivalent to $\omega(q) = \omega(p)$ on $\ker(\alpha_e)$.
If $\omega$ and $\tau$ are cohomology classes, then $\omega+\tau$ and $\omega \tau$ are also cohomology
classes.

\begin{definition}\label{def:cohring}
 The \emph{cohomology ring of $(\Gamma, \alpha)$}, denoted by $H_{\alpha}^*(\Gamma)$, is the subring
 of $\text{Maps}(V,\SS(\ft^*))$ consisting of all the cohomology classes.
\end{definition}

Moreover, $H_{\alpha}^*(\Gamma)$ is a graded ring, with the grading induced by $\SS(\ft^*)$. We say
that $\omega \in H_{\alpha}^*(\Gamma)$ is a class of degree
$2k$ if for every $p \in V$, the polynomial $\omega(p) \in \SS^k(\ft^*)$
is homogeneous of degree $k$. (The fact that the class degree is
twice the polynomial degree is a consequence of the convention that
elements of $\ft^*$ have degree 2.) If $H_{\alpha}^k(\Gamma)$ is the space of
classes of degree $2k$, then
$$H_{\alpha}^*(\Gamma) = \bigoplus_{k \geqslant 0} H_{\alpha}^{2k}(\Gamma) \; .$$

If $\omega \in H_{\alpha}^*(\Gamma)$ and $h \in \SS(\ft^*)$, then $h\omega \in H_{\alpha}^*(\Gamma)$,
hence $H_{\alpha}^*(\Gamma)$ is an $\SS(\ft^*)-$module; it is in fact a graded $\SS(\ft^*)-$module.

\begin{remark}
  The main motivation behind these constructions is the fact that if $M$ is a GKM manifold and
  $\Gamma = \Gamma_M$ is its GKM graph, then $H_T^{odd}(M) = 0$ and $H_T^{2k}(M) \simeq  H_{\alpha}^{2k}(\Gamma)$.
\end{remark}

Let $(\Gamma_0,\alpha)$ be a GKM subgraph of $(\Gamma, \alpha)$. If $f \colon V \to \SS(\ft^*)$
is a cohomology class on $\Gamma$, then the restriction of $f$ to $V_0$ is a cohomology class on $\Gamma_0$.
Therefore the inclusion $i \colon (\Gamma_0, \alpha) \to (\Gamma,\alpha)$ induces a ring morphism
$i^* \colon H_{\alpha}^*(\Gamma) \to H_{\alpha}^*(\Gamma_0)$.

If $\rho = (\varphi,\psi) \colon (\Gamma_1,\alpha_1) \to (\Gamma_2,\alpha_2)$ is an isomorphism of GKM graphs,
define $\rho^* \colon \text{Maps}(V_2,\SS(\ft^*)) \to \text{Maps}(V_1,\SS(\ft^*))$ by
$$(\rho^*(f))(p) = \psi^{-1} (f(\varphi(p))) \; ,$$
for $p \in V_1$, where $\psi^{-1} \colon \SS(\ft^*) \to \SS(\ft^*)$ is the algebra isomorphism extending
the linear isomorphism $\psi^{-1} \colon \ft^* \to \ft^*$. Then $\rho^*$ is a ring isomorphism and
$(\rho^*)^{-1}~=~(\rho^{-1})^*$, but, unless $\psi \colon \ft^* \to \ft^*$ is the identity, $\rho^*$ is
not an isomorphism of $\SS(\ft^*)-$modules.

\subsection{Cohomology of GKM Fiber Bundles}
\label{ssec:3.2}

Let $\pi : (\Gamma, \alpha) \to (B, \alpha_B)$ be a GKM fiber bundle. For a cohomology class
$f \colon V_B \to \SS(\ft^*)$ on the base $(B, \alpha_B)$, define the pull-back
$\pi^*(f) \colon V_{\Gamma} \to \SS(\ft^*)$ by $\pi^*(f)(q) = f(\pi(q))$. Then $\pi^*(f)$ is a
cohomology class on $(\Gamma, \alpha)$, and $\pi$ defines an injective morphism of rings
$\pi^* \colon H_{\alpha_B}^*(B) \to H_{\alpha}^*(\Gamma)$. In particular, $H_{\alpha}^*(\Gamma)$
is an $H_{\alpha_B}^*(B)-$module.

\begin{definition}
  A cohomology class $h \in H_{\alpha}^*(\Gamma)$ is called \emph{basic} if $h \in \pi^*(H_{\alpha_B}^*(B))$.
\end{definition}

Let $(H_{\alpha}^*(\Gamma))_{bas} = \pi^*(H_{\alpha_B}^*(B)) \subseteq H_{\alpha}^*(\Gamma)$.
Then $(H_{\alpha}^*(\Gamma))_{bas}$ is a subring of $H_{\alpha}^*(\Gamma)$, and is isomorphic to
$H_{\alpha_B}^*(B)$. We will identify $H_{\alpha_B}^*(B)$ and $(H_{\alpha}^*(\Gamma))_{bas}$ and
regard $H_{\alpha_B}^*(B)$ as a subring of $H_{\alpha}^*(\Gamma)$.

The next theorem is one of the main results of this paper, and shows
how the cohomology of the total space $\Gamma$ is determined by the
cohomology of the base $B$ and special sets of cohomology classes with certain
properties on fibers.

\begin{theorem}\label{thm:cohtotspace}
  Let $\pi \colon (\Gamma, \alpha) \to (B, \alpha_B)$ be a GKM fiber bundle, and let $c_1,\ldots,c_m$
  be cohomology classes on $\Gamma$ such that, for every $p \in B$, the restrictions of these classes
  to the fiber $\Gamma_{\!p} = \pi^{-1}(p)$ form a basis for the cohomology of the fiber. Then, as
  $H_{\alpha_B}^*(B)-$modules, $H_{\alpha}^*(\Gamma)$ is isomorphic to the free
  $H_{\alpha_B}^*(B)-$module  on $c_1, \ldots, c_m$.
\end{theorem}

\begin{proof}
A linear combination of $c_1,\ldots,c_m$ with coefficients in
$(H_{\alpha}^*(\Gamma))_{bas} \simeq H_{\alpha_B}^*(B)$ is clearly a cohomology class on $\Gamma$.
If such a combination is the zero class, then
$$\sum_{k=1}^m \beta_k(p) c_k(p\,') = 0$$
for every $p \in B$ and $p\,' \in \Gamma_{\!p}$. Since the restrictions of $c_1,\ldots,c_m$ to $\Gamma_{\!p}$
are independent, it follows that $\beta_k(p)=0$ for every $k=1,\ldots,m$. This is valid for all $p\in B$,
hence the classes $\beta_1,\ldots,\beta_m$ are zero. Therefore $c_1, \ldots, c_m$ are independent over
$H_{\alpha_B}^*(B)$, and the free $H_{\alpha_B}^*(B)-$module they generate is a submodule of $H_{\alpha}^*(\Gamma)$.

We prove that this submodule is the entire $H_{\alpha}^*(\Gamma)$. Let $c \in H_{\alpha}^*(\Gamma)$ be a
cohomology class on $\Gamma$. For $p \in B$, the restriction of $c$ to the fiber $\Gamma_{\!p}$ is a
cohomology class on $\Gamma_{\!p}$. Since the restrictions of $c_1,\ldots,c_m$ to $\Gamma_{\!p}$ generate
the cohomology of $\Gamma_{\!p}$, there exist polynomials $\beta_1(p), \ldots, \beta_m(p)$ in $\SS(\ft^*)$
such that, for every $p\,' \in \Gamma_{\!p}$\,,
$$c(p\,') = \sum_{k=1}^m \beta_k(p) c_k(p\,')\; . $$
We will show that the maps $\beta_k \colon B \to \SS(\ft^*)$ are in fact cohomology classes on $B$.
Let $e=p\to q$ be an edge of $B$, with weight $\alpha_e= \alpha_{pq} \in \ft^*$. Let $p\,' \in \Gamma_{\!p}$,
and $q' \in \Gamma_q$ such that $p\,' \to q'$ is the lift of $p \to q$.
Then $\alpha(p\,',q') = \alpha(p,q) = \alpha_e$ and
\begin{align*}
c(q') - c(p\,') = & \sum_{k=1}^m (\beta_k(q) c_k(q') - \beta_k(p) c_k(p\,')) = \\
 = & \sum_{k=1}^m (\beta_k(q)-\beta_k(p)) c_k(p\,') +
 \sum_{k=1}^m \beta_k(q) (c_k(q') - c_k(p\,')) \; .
\end{align*}
Since $c, c_1, \ldots ,c_m$ are classes on $\Gamma$, the differences
$c(q') - c(p\,')$, $c_k(q') - c_k(p\,')$ are multiples of $\alpha_e$, for all $k=1,\ldots, m$.
Therefore, for all $p\,' \in \Gamma_{\!p}$,
$$ \sum_{k=1}^m (\beta_k(q)-\beta_k(p)) c_k(p\,') = \alpha_e \eta(p\,') \; ,$$
where $\eta(p\,') \in \SS(\ft^*)$. We will show that $\eta \colon \Gamma_{\!p} \to \SS(\ft^*)$
is a cohomology class on $\Gamma_{\!p}$.

If $p\,'$ and $p\,''$ are vertices in $\Gamma_{\!p}$, joined by an edge $(p\,',p\,'')$, then
$$ \sum_{k=1}^m (\beta_k(q)-\beta_k(p)) (c_k(p\,'')-c_k(p\,')) =
\alpha_e (\eta(p\,'')- \eta(p\,')) \; .$$
Each $c_k$ is a cohomology class on $\Gamma$, so $c_k(p\,'')-c_k(p\,')$ is a multiple of
$\alpha(p\,',p\,'')$, for all $k=1,\ldots, m$. Then $\alpha_e (\eta(p\,'')- \eta(p\,'))$ is also a
multiple of $\alpha(p\,',p\,'')$. But $\alpha_e =\alpha(p\,',q\,')$ and $\alpha(p\,',p\,'')$ point in
different directions as vectors, so, as linear polynomials, they are relatively prime. Therefore
$\eta(p\,'')- \eta(p\,')$ must be a multiple of $\alpha(p\,',p\,'')$. Therefore $\eta$ is a cohomology
class on $\Gamma_{\!p}$.

The restrictions of $c_1,\ldots,c_m$ form a basis for the cohomology ring of $\Gamma_{\!p}$\,,
hence there exist polynomials $Q_1,\ldots, Q_m \in \SS(\ft^*)$ such that
$$\eta(p\,') = \sum_{k=1}^m Q_k c_k(p\,') \; $$
for all $p\,'\in \Gamma_p$. Then
$$\sum_{k=1}^m (\beta_k(q) - \beta_k(p) - Q_k\alpha_e) c_k = 0 $$
on the fiber $\Gamma_{\!p}$. Since the classes $c_1,\ldots, c_m$ restrict to linearly independent classes
on fibers, it follows that
$$\beta_k(q) - \beta_k(p) = Q_k\alpha_e  \;,$$
hence $\beta_k \in H_{\alpha_B}^*(B)$. Therefore every cohomology class on $\Gamma$ can be written
as a linear combination of classes $c_1,\ldots,c_m$, with coefficients in $H_{\alpha_B}^*(B)$.
\end{proof}

\subsection{Invariant Classes}
\label{ssec:invariant_classes}
In this section we describe a method of constructing global classes $c_1,\ldots,c_m$ with the properties
required by Theorem~\ref{thm:cohtotspace}.

Let $\pi \colon (\Gamma, \alpha) \to (B, \alpha_B)$ be a GKM fiber bundle, with typical fiber
$(F, \alpha_F)$. Let $p$ be a fixed vertex of $B$ and let
$\rho_p = (\varphi_p,\psi_p) \colon (F,\alpha_F) \to (\Gamma_{\!p},\alpha_p)$ be a GKM isomorphism
from $F$ to the fiber above $p$. For a loop $\gamma \in \Omega(p)$, let
$\rho_{\gamma} = (\varphi_{\gamma}, \psi_{\gamma})$ be the GKM automorphism of $(F, \alpha_F)$ determined
by $\gamma$. Let $K = \text{Hol}(F,p)$ be the holonomy subgroup of $\text{Aut}(F,\alpha_F)$ generated by
all automorphisms $\rho_{\gamma}$, and let $f \in (H_{\alpha_F}^*(F))^{K}$ be a cohomology class on
the fiber, invariant under all the automorphisms in $K$. Then
$f_p = (\rho_p^{-1})^*(f) \in H_{\alpha}^*(\Gamma_{\!p})$ is a class on the fiber over $p$, invariant
under all the automorphisms in $\text{Hol}_{\pi}(\Gamma_{\!p}) \subset \text{Aut}(\Gamma_{\!p},\alpha)$.
For any vertex $q \in \Gamma_{\!p}$ we have $f_p(q) \in \SS(\fv_p^*) \subset \SS(\ft^*)$,
where $\fv^*_p$ is the subspace of $\ft^*$ generated by the values of $\alpha$ on the edges
of $\Gamma_{\!p}$.

We will extend the class $f_p$ from the fiber $\Gamma_{\!p}$ to the total space $\Gamma$. Let $p\,'$
be a vertex of $B$, and $\gamma$ a path in $B$ from $p\,'$ to $p$. Let
$\Upsilon_{\gamma}^* \colon H_{\alpha}^*(\Gamma_{\!p}) \to H_{\alpha}^*(\Gamma_{\!p\,'})$ be the
ring isomorphism induced by the GKM graph isomorphism
$\Upsilon_{\gamma} \colon (\Gamma_{\!p\,'}, \alpha) \to (\Gamma_{\!p},\alpha)$. Since $f_p$ is
$\text{Hol}_{\pi}(\Gamma_{\!p})-$invariant, it follows that if $\gamma_1$ and $\gamma_2$ are
two paths in $B$ from $p\,'$ to $p$, then $\Upsilon_{\gamma_1}^*(f_p) = \Upsilon_{\gamma_2}^*(f_{p})$.
We define $f_{p'}=\Upsilon_{\gamma}^*(f_p) \in H_{\alpha}^*(\Gamma_{\!p\,'})$, where $\gamma$ is any
path in $B$ from $p\,'$ to $p$. Then $f_{p'}(q') \in \SS(\ft_{p'}^*) \subset \SS(\ft^*)$ for
every $q' \in \Gamma_{\!p\,'}$.

\begin{proposition}
  Let $c = c_{f,p} \colon V_{\Gamma} \to \SS(\ft^*)$ be defined by $\left.c\right|_{\Gamma_q}= f_q$
  for all $q \in B$. Then $c \in H_{\alpha}^*(\Gamma)$.
\end{proposition}

\begin{proof}
Since the restrictions of $c$ to fibers are classes on fibers, it suffices to show that $c$
satisfies the compatibility conditions along horizontal edges.

Let $(q_1,q_2)$ be a horizontal edge of $\Gamma$ and let $e=(p_1,p_2)$ be the corresponding
edge of $B$. Then
$$c(q_2) - c(q_1) = f_{p_2}(q_2) - f_{p_1}(q_1) = \Psi_e(f_{p_1}(q_1)) - f_{p_1}(q_1)$$
is a multiple of $\alpha_e = \alpha(q_1,q_2)$, because $\Psi_e(x) = x+c(x)\alpha_e$ on $\fv_{p_1}^*$.
\end{proof}

Note that $c$ depends not only on the class $f$ on the typical fiber $F$, but also on the
point $p$ where we start realizing $f$ on $\Gamma$. The choice of $p$ is limited by the fact
that $f$ has to be invariant under the subgroup $\text{Hol}(F,p)$ determined by $p$.

\begin{remark}\label{rem:basis}
  Suppose that the $\SS(\ft^*)-$module $H_{\alpha_F}^*(F)$ has a basis $\{f_1,\ldots,f_m\}$,
  consisting of  $\text{Hol}(F,p)-$invariant classes, for some $p \in B$. Let $c_j = c_{f_j,p}$,
  for $j=1,\ldots,m$. Then the classes $c_1, \ldots,c_m$ have the property that their restrictions
  to each fiber form a basis for the cohomology of the fiber.
\end{remark}

\section{Flag Manifolds as GKM Fiber Bundles}
\label{sec:flags_as_GKM}

Let $G$ be a connected semisimple complex Lie group, let $P$ be a parabolic subgroup of $G$,
and let $M=G/P$ be the corresponding flag manifold. Let $T$ be a maximal compact torus of $G$,
acting on $M$ by left multiplication on $G$. Then $M$ is a GKM space and the
equivariant cohomology ring $H_T^*(M)$ can be computed from the associated GKM graph.

The goal of this section is to briefly review flag manifolds and their GKM graphs. In the last
subsection we will describe the discrete analog of the natural fiber bundle $G/P_1 \to G/P_2$,
with $T \subset P_1 \subset P_2 \subset G$.

\subsection{Flag Manifolds}
\label{ssec:4.1}

In this subsection we review facts about semisimple Lie algebras and flag manifolds. Details
and proofs can be found, for example, in \cite{FH} or \cite{Hu}.

Let $\fg$ be a complex semisimple Lie algebra, $\fh \subset \fg$ a Cartan subalgebra, and $\ft \subset \fh$
a compact real form. Let
$$\fg = \fh \oplus \bigoplus_{\alpha \in \Delta} \fg_{\alpha}$$
be the Cartan decomposition of $\fg$, where $\Delta \subset \ft^*$ is the set of roots.
Let $\Delta^+$ be a choice of positive roots and $\Delta_0 = \{ \alpha_1,\ldots,\alpha_n\} \subset \Delta$
be the corresponding simple roots. The choice of $\Delta^+$ is equivalent to a choice of a
Borel subalgebra $\fb$ of $\fg$,
$$\fb = \fh \oplus \bigoplus_{\alpha \in \Delta^+} \fg_{\alpha} \; .$$
If $G$ is a connected Lie group with Lie algebra $\fg$ and $B$ is the Borel subgroup with
Lie algebra $\fb$, then $M=G/B$ is the manifold of (generalized) complete flags corresponding to $G$.

For a subset $\Sigma \subset \Delta_0$ of simple roots, let
$\langle \Sigma \rangle \subset \Delta^+$ be the set of positive roots that can be written as
linear combinations of roots in $\Sigma$. Then
$$\fp(\Sigma) = \fb \oplus \bigoplus_{\alpha \in \langle \Sigma \rangle} \fg_{-\alpha} \;
= \fh \oplus \bigoplus_{\alpha \in \langle \Sigma \rangle} (\fg_{\alpha} \oplus \fg_{-\alpha})
\oplus \bigoplus_{\alpha \in \Delta^+ \setminus \langle \Sigma \rangle} \fg_{\alpha}
$$
is a Lie subalgebra of $\fg$, and the corresponding Lie subgroup $P(\Sigma) \leqslant G$ is a
parabolic subgroup of $G$. Up to conjugacy, every parabolic subgroup of $G$ is of this form.
The Borel subgroup $B$ corresponds to $\Sigma = \emptyset$, and the whole group $G$ to $\Sigma = \Delta_0$.
The homogeneous space $M=G/P(\Sigma)$ is the manifold of (generalized, partial) flags corresponding
to $G$ and $\Sigma$.

The examples considered in Section~\ref{ssec:motivation} correspond to $G=SL(n,\CC)$.

\subsection{GKM Graphs of Flag Manifolds}

In this subsection we outline the construction of the GKM graph $(\Gamma, \alpha)$
for quotients of parabolic subgroups; more details are available in \cite{GHZ}.

\subsubsection{Weyl groups} For flag manifolds, the construction of the GKM graph
involves Weyl groups and their actions on roots, and we start with a few useful results.
Let $W$ be the Weyl group of $\fg$, generated by reflections $s_\alpha \colon \ft^* \to \ft^*$
for $\alpha \in \Delta_0$.
As a general convention, we will use Greek letters $\alpha$, $\beta$ for roots and axial
functions (whose values are, in this case, roots, and it will be clear from the context
whether $\alpha$ is a root or an axial function), and Roman letters $u$, $v$, $w$, for
elements of the Weyl group $W$. Then $w \beta$ is the element of $\ft^*$ obtained by
applying $w \in W$ to $\beta \in \ft^*$, and $ws_{\beta}$ is the element of the Weyl
group obtained by multiplying $w \in W$ with the reflection $s_{\beta} \in W$ corresponding
to the root $\beta$. Then
%$
$ws_{\beta} = s_{w\beta}w \; ,$
%$
hence two elements of $W$ that differ by a reflection to the left also differ by a
reflection to the right.

For a subset $\Sigma \subset \Delta_0$, let $W(\Sigma)$ be the subgroup of $W$ generated
by reflections $s_{\alpha}$ corresponding to roots $\alpha \in \Sigma$. Then, for a root
$\alpha \in \Delta$, the reflection $s_{\alpha}\in W$ is in $W(\Sigma)$ if and only if
$\alpha \in \langle \Sigma \rangle$ (\cite[1.14]{Hu}).
For subsets $\Sigma_1 \subset \Sigma_2 \subset \Delta_0$, let $W_1 = W(\Sigma_1)$ and
$W_{2}=W(\Sigma_2)$; then $W_1 \leqslant W_2 \leqslant W$.

 \begin{lemma}
 \label{lem:invariant}
   The set $\langle \Sigma_2 \rangle \setminus \langle \Sigma_1 \rangle$ is $W_1-$invariant.
 \end{lemma}

 \begin{proof}
   If $\beta \in \langle \Sigma_2 \rangle \setminus \langle \Sigma_1 \rangle$, then the
positive root $\beta$ is a linear combination of simple roots in $\Sigma_2$, with all
coefficients non-negative. Since $\beta$ is not in $\langle \Sigma_1 \rangle$, there exists
at least one simple root, say $\alpha_i$, that is not in $\Sigma_1$ and appears in $\beta$
with a strictly positive coefficient. If $\alpha \in \Sigma_1$, then
$s_{\alpha}\beta = \beta - n_{\beta,\alpha}\alpha$, with $n_{\beta,\alpha} \in \ZZ$.
Then $s_{\alpha}\beta$ and $\beta$ have the same coefficients in front of the simple
roots not in $\Sigma_1$. In particular, $\alpha_i$ appears in $s_{\alpha}\beta$ with a
strictly positive coefficient, which proves that $s_{\alpha}\beta$ is a positive root.
The simple roots appearing in $\alpha$ and $\beta$ are all in $\Sigma_2$, hence
$s_{\alpha}\beta \in \langle \Sigma_2 \rangle$, and as $\alpha_i$ is not in $\Sigma_1$,
it follows that $s_{\alpha}\beta \in \langle \Sigma_2 \rangle \setminus \langle \Sigma_1 \rangle$.
Since $W_1$ is generated by the reflections $s_{\alpha}$ with $\alpha \in \Sigma_1$, we conclude
that $\langle \Sigma_2 \rangle \setminus \langle \Sigma_1 \rangle$ is $W_1-$invariant.
 \end{proof}

Let $w \in W_2$ and let $w=s_{\beta_1} \dotsb s_{\beta_m}$ be a decomposition of $w$ into
simple reflections, with $\beta_i \in \Sigma_2$ for all $i=1,...,m$. If
$\alpha \in \langle \Sigma_1 \rangle$ and $\beta \in \langle \Sigma_2 \rangle
\setminus \langle \Sigma_1 \rangle$ then
$$s_{\beta} s_{\alpha} = s_{\alpha} s_{s_{\alpha}\beta} \; , $$
and $s_{\alpha}\beta \in \langle \Sigma_2 \rangle \setminus \langle \Sigma_1 \rangle$.
We can therefore push all the reflections coming from roots in $\langle \Sigma_1 \rangle$
to the left, and get $w=us_{\beta_1'}\dots s_{\beta_k'}$ with $u \in W_1$ and
$\beta_1',\ldots,\beta_k' \in \langle \Sigma_2 \rangle \setminus \langle \Sigma_1 \rangle$.
We can also push all the reflections coming from roots in $\langle \Sigma_1 \rangle$ to the
right, and get $w=s_{\beta_1''}\dots s_{\beta_k''}u$ with $u \in W_1$ and
$\beta_1'',\ldots,\beta_k'' \in \langle \Sigma_2 \rangle \setminus \langle \Sigma_1 \rangle$.

\subsubsection{Quotients of parabolic subgroups}

Let $\Sigma_1 \subset \Sigma_2 \subset \Delta_0$ be subsets of simple roots and
$B \leqslant P(\Sigma_1):=P_1 \leqslant P(\Sigma_2):=P_2 \leqslant G$ the corresponding
parabolic subgroups. The compact torus $T$ with Lie algebra $\ft$ acts on $M=P_2/P_1$ by
left multiplication on $P_2$, and the space $M=P_2/P_1$ is a GKM space,
isomorphic to $G'/P'$ for a Levi subgroup $G'$ of $P_1$. All flag manifolds
are of this type, corresponding to $\Sigma_2 = \Delta_0$.

We describe now the GKM graph $(\Gamma, \alpha)$ associated to $M=P_2/P_1$. The fixed point
set $M^T$ is identified with the set of right cosets
$$W_2/W_1 = \{ vW_1 \; | \; v \in W_2\; \}  = \{ [v] \; | \; v \in W_2\; \}\; ,$$
where $[v] = vW_1$ is the right $W_1-$coset containing $v \in W_2$. Vertices $[w],[v]$
are joined by an edge if and only if $[v] = [ws_{\beta}]$ for some
$\beta \in \langle \Sigma_2 \rangle \setminus \langle \Sigma_1 \rangle$.
If $[w\sigma_{\beta}]=[w]$, then $\sigma_{\beta} \in W_1$, which is impossible if
$\beta \in \langle \Sigma_2 \rangle \setminus \langle \Sigma_1 \rangle$,
because the only reflections in $W_1$ are those associated to roots in $\Sigma_1$.
Therefore the endpoints of an edge are distinct and the graph has no loops.
For $w \in W_2$ and $\beta \in \langle \Sigma_2 \rangle \setminus \langle \Sigma_1 \rangle$, the
edge $e = ([w] \to [ws_{\beta}] = [s_{w\beta}w])$ is labeled by
$\alpha_e = \alpha([w],[ws_{\beta}]) = w\beta$.

We show that the label $\alpha_e$ is independent of the representative $w\in W_2$:
if $[w']=[w]$ and $[ws_{\beta}] = [w's_{\gamma}]$ with
$\beta, \gamma \in \langle \Sigma_2 \rangle \setminus \langle \Sigma_1 \rangle$,
then there exist $w_1, w_2 \in W_1$ such that $w'=ww_1$ and $w's_{\gamma} = ws_{\beta}w_2$.
Then $s_{\beta}s_{w_1\gamma} = w_2w_1^{-1} \in W_1$, which implies $w_1\gamma = \pm \beta$.
Since $\langle \Sigma_2 \rangle \setminus \langle \Sigma_1 \rangle$ is $W_1-$invariant, it
follows that $w_1\gamma = \beta$ and therefore $w'\gamma = ww_1\gamma = w\beta$.

The connection along the edge $e = ([w], [ws_{\beta}])$ sends the edge $e' = ([w],[ws_{\beta'}])$
to the edge $e'' = ([ws_{\beta}], [ws_{\beta}s_{\beta'}])$.

Then $(\Gamma(W_2/W_1), \alpha)$ is the GKM graph of the GKM space $M=P_2/P_1$. We will refer
to it simply as $W_2/W_1$, and it will be clear from the context when we mean the GKM graph,
when just the graph, and when just the vertices.

\begin{example}
\label{exm:particular_cases}
 We describe the particular cases when $P_2=G$ or $P_1=B$, or both.

  For $M=G/B$ we have $\Sigma_1=\emptyset$, $\Sigma_2 = \Delta_0$, $W_1=\{1\}$ and $W_2=W$,
hence $W_2/W_1= W$. Vertices $w,v \in W$ of the corresponding GKM graph $\Gamma(W)$ are joined
by an edge if and only if $w^{-1}v = s_{\beta}$ for some $\beta \in \Delta^+$ (or,
equivalently, if $v = ws_{\beta} = s_{w\beta}w$), and the edge $w \to ws_{\beta} = s_{w\beta}w$
is labeled by $w\beta$.

    For $M=P(\Sigma)/B$, we have $\Sigma_1= \emptyset$, $\Sigma_2 = \Sigma \subset \Delta_0$,
$W_2=W(\Sigma)$, and $W_1 = \{ 1 \}$. The GKM graph $\Gamma(W(\Sigma))$ is the induced subgraph
of $\Gamma(W)$ with vertex set $W(\Sigma)$: vertices $w,v \in W(\Sigma)$ are joined by an edge
in $\Gamma(W(\Sigma))$ if and only if they are joined by an edge in $\Gamma(W)$. That happens
if $v=ws_{\beta} = s_{w\beta}w$ for some $\beta \in \langle \Sigma \rangle$. The edge
$w \to ws_{\beta} = s_{w\beta}w$ is labeled by $w\beta$.

  For $M=G/P(\Sigma)$, we have $\Sigma_2 = \Delta_0$ and $\Sigma_1= \Sigma \subset \Delta_0$.
The GKM graph is a graph with vertex set $W/W(\Sigma)$. Vertices $[w],[v] \in W/W(\Sigma)$
are joined by an edge if and only if $w^{-1}v = s_{\beta}$ for some
$\beta \in \Delta^+ \setminus \langle \Sigma \rangle$; equivalently, if
$v=ws_{\beta} = s_{w\beta}w$. The edge $w \to ws_{\beta} = s_{w\beta}w$ is labeled by $w\beta$.
\end{example}

\subsection{GKM Fiber Bundles of Flag Manifolds}
\label{ssec:4.2}

Let $\Sigma_1 \varsubsetneq \Sigma_2 \subset \Delta_0$ be, as above, subsets of simple roots,
and let $W_1 = W(\Sigma_1)$ and $W_2 = W(\Sigma_2)$ be the corresponding subgroups of $W$.
For an element $w\in W$, let $wW_1$ be its class in $W/W_1$, and $wW_2$ its class in $W/W_2$.
One has a natural map $\pi \colon W/W_1 \to W/W_2$, given by $\pi(wW_1) = wW_2$, from the
vertices of $\Gamma(W/W_1)$ to the vertices of $\Gamma(W/W_2)$. If $\Sigma_2 = \Delta_0$,
then the base $W/W_2$ is just a point and the map $\pi$ is trivial. For the rest of this section
we will assume that $\Sigma_2 \varsubsetneq \Delta_0$. The goal of this section is to show that
$\pi$ is a GKM fiber bundle between the corresponding GKM graphs.

\begin{theorem}
  The projection $\pi \colon W/W_1 \to W/W_2$ is a GKM fiber bundle with typical fiber $W_2/W_1$.
\end{theorem}

\begin{proof}
Let $wW_1$ be a vertex of $W/W_1$ and let $e=(wW_1,ws_{\beta}W_1)$ be an edge of $W/W_1$,
with $\beta \in \Delta^+ \setminus \langle \Sigma_1 \rangle$. This edge is vertical if and only
if $s_{\beta} \in W_2$, and this happens exactly when $\beta \in \langle \Sigma_2 \rangle$.
Therefore the vertical edges at $wW_1$ are
$$E_{wW_1}^{\bot} = \{ (wW_1,ws_{\beta}W_1) \; | \; \beta \in \langle \Sigma_2 \rangle \setminus
\langle \Sigma_1 \rangle \; \} \; ,$$
and the horizontal edges are
$$H_{wW_1} = \{ (wW_1,ws_{\beta}W_1) \; | \; \beta \in \Delta^+ \setminus \langle \Sigma_2 \rangle
\; \} \; .$$
If $(wW_1,ws_{\beta}W_1)$ is a horizontal edge, then $(wW_2,ws_{\beta}W_2)$ is an edge of $W/W_2$,
hence $\pi$ is a morphism of graphs, and $(d\pi)_{wW_1} \colon H_{wW_1} \to E_{wW_2}$, is defined by
$$(d\pi)_w(wW_1,ws_{\beta}W_1) = (wW_2,ws_{\beta}W_2)\; .$$
It is clear that $(d\pi)_{wW_1}$ is a bijection, hence $\pi$ is a fibration of graphs.

Next we show that $\pi$ is a GKM fibration. Let $e=(wW_2,ws_{\beta}W_2)$ be an edge of $W/W_2$,
with $\beta \in \Delta^+ \setminus \langle \Sigma_2 \rangle$. If $vW_1$ is a vertex of $W/W_1$ in
the fiber above $wW_2$, then $v=wu$, for some $u \in W_2$. Let $\beta'= u^{-1}\beta$.
By Lemma~\ref{lem:invariant} applied to the pair $(\Delta_0,\Sigma_2)$ corresponding to $(W,W_2)$,
the set $\Delta^+ \setminus \langle \Sigma_2\rangle$ is $W_2-$invariant, hence
$\beta' \in \Delta^+ \setminus \langle \Sigma_2 \rangle$. Therefore $\tilde{e}=(vW_1,vs_{\beta'}W_1)$
is an edge of $W/W_1$. Since
$$\pi(vs_{\beta'}W_1) = vs_{\beta'}W_2 = wus_{u^{-1}\beta}W_2 = ws_{\beta}uW_2 = ws_{\beta}W_2\; , $$
it follows that $\tilde{e}$ is the lift of $e$ at $vW_1$. Moreover, if $\alpha_1$ and $\alpha_2$ are
the axial functions on $W/W_1$ and $W/W_2$, respectively, then
$$\alpha_1(vW_1,vs_{\beta'}W_1) = v\beta' = wuu^{-1}\beta = w\beta = \alpha_2(wW_2,ws_{\beta}W_2)\; ,$$
hence the axial functions are compatible with $\pi$.

Let $e = (vW_1,vs_{\beta}W_1)$ and $e' = (vW_1,vs_{\beta'}W_1)$ be edges of $W/W_1$. The connection
$\nabla_1$ along $e$ moves $e'$ to $e'' = (vs_{\beta}W_1,vs_{\beta}s_{\beta'}W_1)$. If
$\beta' \in \Delta^+ \setminus \langle \Sigma_2\rangle$, then both $e'$ and $e''$ are horizontal,
and if $\beta' \in  \langle \Sigma_2\rangle \setminus \langle \Sigma_1 \rangle$, then both are vertical.
Hence the connection along any edge of $W/W_1$ moves horizontal edges to horizontal edges and vertical
edges to vertical edges. Moreover, if both $e$ and $e'$ are horizontal (and hence so is $e''$), then
the connection $\nabla_2$ along the projection of $e$ moves the projection of $e'$ to the projection
of $e''$, which shows that the restriction of $\nabla_1$ to horizontal edges is compatible with
$\nabla_2$, and we have shown that $\pi$ is a GKM fibration.

Finally, we prove that $\pi$ is a GKM fiber bundle. Let $p=wW_2$ and $q=ws_{\beta}W_2$ be two
adjacent vertices of $W/W_2$, with $\beta \in \Delta^+ \setminus \langle \Sigma_2 \rangle$.
A straightforward computation shows that the transition map $\Phi_{p,q} \colon \pi^{-1}(p) \to \pi^{-1}(q)$
is given by
$$\Phi_{p,q} (vW_1) = s_{w\beta}vW_1 \; ,$$
and hence, if $e' = (vW_1,vs_{\beta'}W_1)$ is an edge of $\pi^{-1}(p)$, then
$$e'' = (\Phi_{p,q}(vW_1), \Phi_{p,q}(vs_{\beta'}W_1)) = (s_{w\beta}vW_1, s_{w\beta}vs_{\beta'}W_1)$$
is an edge of $\pi^{-1}(q)$. Therefore $\Phi_{p,q}$ is a morphism of graphs, hence an isomorphism,
with inverse $\Phi_{p,q}^{-1} = \Phi_{q,p}$. In addition, the connection $\nabla_1$ along the lift
of $e=(p,q)$ at $vW_1$ moves $e'$ to $e''$. Moreover
$$\alpha_1(e'') = s_{w\beta}v\beta' = s_{w\beta}(\alpha_1(e'))\, ,$$
hence, if $\Psi_{p,q} \colon \ft^* \to \ft^*$ is given by $\Psi_{p,q}(x) = s_{w\beta}(x)$, then
its induced restriction and co-restriction $\Psi_{p,q} \colon \fv_p^* \to \fv_q^*$ is compatible
with $\Phi_{p,q}$. This proves that
$$(\Phi_{p,q}, \Psi_{p,q}) \colon (W/W_1)_p \to (W/W_1)_q$$
is an isomorphism of GKM graphs, hence the fibers are canonically isomorphic, through an
isomorphism compatible with the connection of $\Gamma_1$. We conclude that $\pi$ is a GKM fiber bundle.

All that remains is to show that the fibers are isomorphic, as GKM graphs, to $W_2/W_1$. Let $p$ be a vertex of $W/W_2$
and $w\in W$ a representative for $p$. Let $\varphi_w \colon W_2/W_1 \to \pi^{-1}(p)$,
$\varphi_w(vW_1) = wvW_1$ and $\psi_w$ the restriction and co-restriction of
$\psi_w \colon \ft^* \to \ft^*$, $\psi_w (x) = wx$. Note that $\varphi_w$ and $\psi_w$
depend not just on the class $p$, but on the \emph{particular} representative $w$. If
$e=(vW_1,vs_{\beta}W_1)$ is an edge of $W_2/W_1$, with $\beta \in \langle \Sigma_2 \rangle \setminus
\langle \Sigma_1 \rangle$, then $e'=(\varphi_w(vW_1), \varphi_w(vs_{\beta}W_1) = (wvW_1,wvs_{\beta}W_1)$
is an edge of the fiber, and
$$\alpha_1(e') = wv\beta = \psi_v(\alpha(e))\; .$$
It is not hard to see that $(\varphi_w,\psi_w) \colon W_2/W_1 \to \pi^{-1}(p)$ is in fact an
isomorphism of GKM graphs, and this concludes the proof of the theorem.
\end{proof}

The example considered in Section~\ref{ssec:sum_example} is the particular case of a root system of
type $A_{n-1}$, with $\Sigma_1 = \emptyset$ and $\Sigma_2 = \Delta_0 \setminus \{ \alpha_1\}$.
The fiber bundle $\mathcal{F}l_4(\CC) \to \mathcal{G}r_2(\CC^4)$ shown in Figure~\ref{fig:j42}
corresponds to the root system $A_3$, with $\Delta_0 = \{ \alpha_1, \alpha_2, \alpha_3\}$,
$\Sigma_1 = \emptyset$ and $\Sigma_2 =\{ \alpha_1, \alpha_3\}$.

\subsection{Holonomy Subgroup}
\label{ssec:4.3}

In this section we determine the holonomy subgroup of $\text{Aut}(W_2/W_1,\alpha)$ determined by
loops in the base $W/W_2$.

Let $w \in W_2$, let $\Phi_w \colon W_2/W_1 \to W_2/W_1$, $\Phi_w(uW_1) = wuW_1$, and
$\Psi_w \colon \ft^* \to \ft^*$, $\Psi_w(\beta) = w\beta$. Then $\Upsilon_w = (\Phi_w,\Psi_w)
\colon W_2/W_1 \to W_2/W_1$ is a GKM automorphism. Moreover, the map $\Upsilon \colon W_2 \to
\text{Aut}(W_2/W_1,\alpha)$, $\Upsilon(w) = \Upsilon_w$ is a morphism of groups with kernel
included in $W_1$. When $W_1$ is a normal subgroup of $W_2$, the kernel \emph{is} $W_1$, and
then the image $\Upsilon(W_2)$ is isomorphic
with the quotient group $W_2/W_1$.

\begin{proposition}
The holonomy subgroup of $\text{Aut}(W_2/W_1,\alpha)$ is $\Upsilon(W_2)$.
\end{proposition}

\begin{proof}
For $v_0 \in W$ let $\pi^{-1}(v_0W_2) \subset W/W_1$ be the fiber through $v_0W_2$, identified
with $W_2/W_1$ by $(\varphi_{v_0}, \psi_{v_0}) \colon W_2/W_1 \to \pi^{-1}(v_0W_2)$.

Let $\gamma \in \Omega(v_0W_2)$ be a loop in $W/W_2$ based at $v_0W_2$, given by
$$v_0W_2 \to v_1W_2 \to \dotsb \to v_{m-1}W_2 \to v_mW_2=v_0W_2\; ,$$
where $v_k = v_{k-1}s_{\beta_k}$, with $\beta_k \in \Delta^{+} \setminus \langle \Sigma_2 \rangle$
for $k=1,\ldots,m$, and let $w=v_0^{-1}v_m$. Then
$w= s_{\beta_1}\dotsb s_{\beta_m}$, and since $\gamma$ is a loop, we have
$w \in W_2$.

Let $\varphi_{\gamma} \colon W_2/W_1 \to W_2/W_1$ be the map
$$\varphi_{\gamma} = \varphi_{v_0}^{-1} \circ \Phi_{\gamma} \circ \varphi_{v_0} = \varphi_{v_0}^{-1}
\circ \Phi_{v_{m-1}W_2,v_mW_2} \circ \dotsb \circ \Phi_{v_{0}W_2,v_1W_2} \circ \varphi_{v_0}\; . $$
Then
$$\Phi_{v_{0}W_2,v_1W_2} \circ \varphi_{v_0}(uW_1) = s_{v_0\beta_1}v_0uW_1 = v_0s_{\beta_1}uW_1 =
\varphi_{v_1}(uW_1) \; .$$
Continuing with the other edges of $\gamma$, we get
$$\varphi_{\gamma}(uW_1) = \varphi_{v_0}^{-1} \varphi_{v_m}(uW_1) = \Phi_w(uW_1)\; ,$$
hence $\varphi_{\gamma} = \Phi_w$. Similarly, $\psi_{\gamma} = \psi_w$, and hence $\rho_{\gamma} = \Upsilon_w$.
We conclude that
$$\text{Hol}(W_2/W_1,v_0W_2) \subset \Upsilon(W_2)\; .$$

We now show that for every $v \in W_2$, there exists a loop $\gamma$ in $W/W_2$, starting and
ending at $v_0W_2$, and such that $\rho_{\gamma}=\Upsilon(v)$.

Let $\alpha_i \in \Sigma_2 \varsubsetneq \Delta_0$. The Weyl group $W$ acts transitively on $\Delta$,
hence there exists $w \in W$ such that $w \alpha_i \in \Delta^{+} \setminus \langle \Sigma_2 \rangle$.
Let $w=uv$  be a decomposition of $w$ such that $u \in W_2$ and $v=s_{\beta_1}\dotsb s_{\beta_m}$
with $\beta_1,\ldots,\beta_m \in \Delta^{+} \setminus \langle \Sigma_2 \rangle$. Then
$u^{-1}w\alpha_i \in \Delta^{+} \setminus \langle \Sigma_2 \rangle$, because
$\Delta^{+} \setminus \langle \Sigma_2 \rangle$ is $W_2-$invariant. Consider the path
$\gamma$ in $W/W_2$ that starts with
$$v_0W_2 \to v_0 s_{\beta_m}W_2 \to \dotsb \to v_0 s_{\beta_m} \dotsb s_{\beta_1}W_2 = v_0v^{-1}W_2 \; , $$
continues with
$$v_0 v^{-1}W_2 \to v_0 v^{-1} s_{u^{-1}w\alpha_i}W_2 \to v_0 v^{-1} s_{u^{-1}w\alpha_i} s_{\beta_1}W_2
\to v_0 v^{-1} s_{u^{-1}w\alpha_i} s_{\beta_1}s_{\beta_2}W_2\; ,$$
and ends with
$$v_0 v^{-1} s_{u^{-1}w\alpha_i} s_{\beta_1}s_{\beta_2}W_2 \to \dotsb \to
v_0 v^{-1} s_{u^{-1}w\alpha_i} s_{\beta_1}s_{\beta_2}\dotsb s_{\beta_m}W_2 =
v_0 v^{-1} s_{u^{-1}w\alpha_i} vW_2 \; .$$
This path is a loop because $v_0 v^{-1} s_{u^{-1}w\alpha_i} v = v_0s_{\alpha_i}$ and
$\alpha_i \in \Sigma_2$, and
$$\rho_{\gamma} = \Upsilon_{v_0v_0^{-1}s_i} = \Upsilon(s_i)\; .$$

Since $W_2$ is generated by $s_i=s_{\alpha_i}$ for $\alpha_i \in \Sigma_2$, we conclude that
$$\text{Hol}(W_2/W_1,v_0W_2) = \Upsilon(W_2)\; ,$$
and the holonomy group of the typical fiber does not depend on the base point.
\end{proof}

\subsection{Bases of Invariant Classes}
\label{ssec:4.1.5}

We use the GKM graph of $M = G/B$ to describe equivariant cohomology classes in $H_T^*(M)$.
The ring $H_{\alpha}^*(W)$ consists of the maps $f \colon W \to \SS(\ft^*)$ such that
$$f(ws_{\beta}) - f(w) \in (w\beta) \SS(\ft^*)$$
for every $w \in W$ and $\beta \in \Delta^+$.

The Weyl group action on $\ft^*$ induces an action of $W$ on $H_{\alpha}^*(W)$,
given by
$$w \cdot f = f^{w} \colon W \to \SS(\ft^*)\; , \quad f^{w}(v) = w^{-1} f(wv) \; .$$

Let $K$ be a compact real form of $G$ containing $T$. Then (see, for example,
\cite[Section 4.7]{GS}) the subring of $W-$invariant classes is
$$H_\alpha^*(W)^W \simeq H_T^*(M)^W \simeq H_K^*(M) = H_K^*(G/B) =H_T^*(K/T) \simeq \SS(\ft^*)\; ,$$
An explicit ring isomorphism from $\SS(\ft^*)$ to $H_{\alpha}^*(W)^W$ is given by
\begin{equation}
 \label{eq:equiv_charact_hom}
c_T \colon \SS(\ft^*) \to H_{\alpha}^*(W)^W, c_T(q)(v) = v\cdot q\; ,
\end{equation}
for all $q \in \SS(\ft^*)$ and $v \in W$. The inverse is
$c_T^{-1} \colon H_{\alpha}^*(W)^W \to \SS(\ft^*)$, $c_T^{-1}(f) = f(1)$.

We will show in Section~\ref{ssec:sym_schubs} that the $\SS(\ft^*)-$module $H_\alpha^*(W)$ has bases
consisting of $W-$invariant classes. The isomorphism $c_T$ establishes an explicit
correspondence between
such bases and $\SS(\ft^*)^W-$module bases of $\SS(\ft^*)$.

\begin{theorem}
\label{prop:bases_over_invariants}
Let $q_1, \ldots, q_N$ be elements of $\SS(\ft^*)$ and $f_i = c_T(q_i)$, $i=1,\ldots,N$
the corresponding $W-$invariant classes.
Then $\{ f_1, \ldots, f_N\}$ is a basis of $H_\alpha^*(W)$ over $\SS(\ft^*)$  if and only
if $\{q_1, \ldots, q_N\}$ is a basis of $\SS(\ft^*)$ over $\SS(\ft^*)^W$.
\end{theorem}

\begin{proof}
  Assume first that $\{ f_1, \ldots, f_N\}$ is a basis of $H_\alpha^*(W)$ over $\SS(\ft^*)$.

 Suppose that $a_1, \ldots, a_N$ are elements of $\SS(\ft^*)^W$ such that
  $$a_1 q_1 + \dotsb + a_N q_N = 0\; .$$
  Then for every $v \in W$ we have
  $$
    v\cdot (a_1 q_1 + \dotsb + a_N q_N)  = 0 \Longrightarrow
    a_1 f_1(v) + \dotsb + a_N f_N(v)  = 0 \;, $$
and since this is valid for every $v \in W$, we conclude that
$$
    a_1 f_1 + \dotsb + a_N f_N  = 0\; .
  $$
  But the classes $f_1, \ldots, f_N$ are independent, hence $a_1= \dotsb = a_N = 0$.
Therefore $q_1,\ldots,q_N$ are linearly independent over $\SS(\ft^*)^W$.

  Let $q \in \SS(\ft^*)$. Then $c_T(q) \in H_{\alpha}^*(W)$, hence there exist
$a_1,\ldots, a_N$ in $\SS(\ft^*)$ such that
  $$c_T(q) = a_1 f_1 + \dotsb + a_N f_N\; .$$
  Then for every $v \in W$ we have
  \begin{align*}
    c_T(q)(v^{-1}) & = a_1 f_1(v^{-1}) + \dotsb + a_N f_N(v^{-1})
\Longrightarrow \\ v^{-1}\cdot q & = a_1 v^{-1}\cdot q_1 + \dotsb + a_N v^{-1} \cdot q_N \Longrightarrow \\
    q & = (v\cdot a_1) \, q_1 + \dotsb + (v \cdot a_N) \, q_N \; .
  \end{align*}

  Averaging over $W$ we get
  $$q = b_1 q_1+ \dotsb +b_N q_N \, ,$$
  where for each $k = 1,\ldots, N$,
  $$b_k = \frac{1}{|W|} \sum_{v \in W} v\cdot a_k$$
  is an element of $\SS(\ft^*)^W$. This proves that $q_1, \ldots, q_N$ also generate $\SS(\ft^*)$ over $\SS(\ft^*)^W$,
  and therefore $\{q_1, \ldots, q_N\}$ is a basis of $\SS(\ft^*)$ over $\SS(\ft^*)^W$.

  Conversely, assume now that $\{q_1, \ldots, q_N\}$ is a basis of $\SS(\ft^*)$ over $\SS(\ft^*)^W$.

  Let $\{\sigma_1, \ldots, \sigma_N\}$ be a basis of $H_{\alpha}^*(W)$ consisting of $W-$invariant classes.
  There must be exactly $N$ such classes, because by the first part $\{r_i=\sigma_i(1)\, | \, i=1,..,N\}$
  is a basis of
  $\SS(\ft^*)$ over $\SS(\ft^*)^W$, and all bases of a free module over a commutative ring have the same
  number of elements.

  Let $A \in GL_N(\SS(\ft^*)^W) \subset GL_N(\SS(\ft^*))$ be the change-of-basis matrix from the
  basis $\{r_1, \ldots, r_N\}$ to the
  basis $\{q_1, \ldots, q_N\}$:
  $$q_i = a_{1i} r_1 + \dotsb + a_{Ni} r_N$$
  for all $i=1\ldots, N$. Since the entries of $A$ are $W-$invariant, for $v \in W$ we have
$$f_i(v)  = v\cdot q_i = a_{1i} v\cdot r_1  + \dotsb + a_{Ni} v\cdot r_N = a_{1i} \sigma_1(v)  + \dotsb + a_{Ni} \sigma_N(v)$$
  and therefore
  $$f_i = a_{1i} \sigma_1  + \dotsb + a_{Ni} \sigma_N$$
  for all $i=1\ldots,N$. Since $\{\sigma_1,\ldots,\sigma_N\}$ is a basis and $A$ is
  invertible, it follows that $\{f_1,\ldots,f_N\}$ is also a basis, and that concludes the proof.
\end{proof}

\section{Fibrations of Classical Groups}
\label{sec:classical}

In this section we consider the GKM bundle $W \to W/W_S$ when  $S=\Delta_0 \setminus \{ \alpha_1\}$,
where $\Delta_0$ is the set of simple roots for a classical root
system and $\alpha_1$ is one of
the endpoint roots in the Dynkin diagram. By recursively applying
Theorem~\ref{thm:cohtotspace}, we construct a basis of $H_{\alpha}^*(W)$
consisting of $W-$invariant classes.

\subsection{Type A}
The set of simple roots of $A_n$ (for $n \geqslant 2$) is $\Delta_0 = \{\alpha_1,\ldots,\alpha_n\}$,
where $\alpha_i = x_i - x_{i+1}$, for $i=1,\ldots,n$. The set of positive roots is
$$\Delta^{+} = \{ x_i - x_j \; | \; 1 \leqslant i < j \leqslant n+1\}$$
and $x_i-x_j = \alpha_i + \ldots + \alpha_{j-1}$. If $S= \{\alpha_2, \ldots,\alpha_n\}$, then
$$\langle S \rangle = \{ x_i -x_j \; | \; 2 \leqslant i < j \leqslant n+1\} \; ,$$
is the set of positive roots for a root system of type $A_{n-1}$, and
$$\Delta^+ \setminus \langle S \rangle = \{ \beta_j  \; | \;
\beta_j = x_1- x_j, 2 \leqslant j \leqslant n+1\} = \{ \alpha_1 + \dotsb + \alpha_j \; | \;
1 \leqslant j \leqslant n\} \; .$$

Let
$$
  \omega_1 =  [id] \quad \text{and} \quad
  \omega_j =  [s_{\beta_j}]\; , \text{ for } 2 \leqslant j \leqslant n+1 \; .
$$

Then $W/W_S = \{\omega_1,\ldots,\omega_{n+1} \}$, and the graph structure of $W/W_S$
is that of a complete graph with $n+1$ vertices. If $\tau \colon W/W_S \to \ft^*$ is
given by $\tau(\omega_i) = x_i$ for all $i=1,\ldots,n+1$, then the axial function
$\alpha$ on $W/W_S$ is given by
$$\alpha(\omega_i, \omega_j) = \tau(\omega_i) - \tau(\omega_j) = x_i - x_j$$
and $\tau \in H_{\alpha}^1(W/W_S)$ is a class of degree 1. Using a Vandermonde
determinant argument, one can show that the classes $\{ 1, \tau, \ldots, \tau^n\}$ are
linearly independent over $\SS(\ft^*)$, and in fact form a basis of the free $\SS(\ft^*)-$module
$H_{\alpha}^*(W/W_S)$.

The Weyl group $W$ is isomorphic to the symmetric group $S_{n+1}$, acting on roots by
$$w \cdot (x_i-x_j) = x_{w(i)} - x_{w(j)}\; .$$
The simple reflection $s_i$ acts as the transposition $(i,i+1)$, and, more generally,
the reflection associated to the root $x_i-x_j$ acts as the transposition $(i,j)$.
The subgroup $W_S$ is the subgroup of $W=S_{n+1}$ consisting of the permutations that
fix the element 1. With the identification $W/W_S \simeq K_{n+1}$, the projection
$\pi \colon W \to W/W_S$ is the map $\pi \colon S_{n+1} \to K_{n+1}$, $\pi(w) = w(1)$.

\begin{remark}
  This is essentially the example discussed in Section ~\ref{ssec:sum_example},
and corresponds to the fiber bundle of complete flags over a projective space. The group $G$ is
$SL_{n+1}(\CC)$,
the Borel subgroup $B$ is the subgroup of upper triangular matrices, and the parabolic subgroup $P$ is
the subgroup of $G$ consisting of block-diagonal matrices, with one block of size $1\times 1$ and a
second block of size $n \times n$. Then $G/B \simeq \mathcal{F}l(\CC^{n+1})$ and $G/P \simeq \CC P^n$.
The projection $\pi \colon \mathcal{F}l(\CC^{n+1}) \to \CC P^n$ sends the flag
    $$V_{\bullet} \colon V_1 \subset V_2 \subset \dotsb \subset V_{n} \subset \CC^{n+1}$$
to $\pi(V_{\bullet}) = V_1$. For an element $L \in \CC P^{n}$, hence a one-dimensional subspace
of $\CC^{n+1}$, the fiber $\pi^{-1}(L)$ is diffeomorphic to
$\mathcal{F}l(\CC^{n+1}/L) \simeq \mathcal{F}l(\CC^{n})$.
\end{remark}

For a multi-index $I=[i_1,\ldots,i_{n}]$ of non-negative integers, we define
$$\textbf{x}^I = {x_1}^{i_1} {x_2}^{i_2}\dotsb x_{n}^{i_{n}}$$
and let $c_I = c_T(\textbf{x}^I)$ be the corresponding $W-$invariant class
$c_I \colon S_{n+1} \to \SS(\ft^*)$,
$$c_I(u) = u \cdot \textbf{x}^I = x_{u(1)}^{i_1} \dotsb x_{u(n)}^{i_{n}}\; ;$$
then $\textbf{x}^I = c_I(id)$, where $id$ is the identity element of the Weyl group $W=S_{n+1}$.
We will construct a basis of the $\SS(\ft^*)-$module $H_{\alpha}^*(W)$ consisting of
classes of the type $c_I$ for specific indices $I$.

Consider the GKM fiber bundle $\pi \colon S_3 \to K_3$, $\pi(u) = u(1)$.
The fiber $\pi^{-1}(3)$ is canonically isomorphic to $S_2$, and since $S_2 \simeq K_2$,
the cohomology of $S_2$ is a free $\SS(\ft^*)-$module with a basis given by the invariant
classes $c_{[0]}$ and $c_{[1]}$.  The invariant class $c_{[0]}$ on this fiber is extended,
using transition maps  between fibers to the constant class $c_{[0,0]} \equiv 1$ on the
total space. The invariant class $c_{[1]}$ extends to the class $c_{[0,1]}$;
the shift in index is due to the fact that the axial functions on fibers are different. The cohomology
of the base $K_3$ is generated, over $\SS(\ft^*)$, by $1$, $\tau$, and $\tau^2$, and these
classes lift to basic classes $c_{[0,0]}$, $c_{[1,0]}$, and $c_{[2,0]}$ on $S_3$.
Theorem~\ref{thm:cohtotspace} implies that the cohomology of $S_3$ is a free $\SS(\ft^*)-$module,
with a basis given by
$$\{c_I \; | \; I = [i_1,i_2], 0 \leqslant i_1 \leqslant 2, 0 \leqslant i_2 \leqslant 1 \} \; .$$
Their values on $W(A_2) =S_3$ are given in Table~\ref{tbl:inv_class_S3}.

\begin{table}[h]
  \begin{tabular}{|c|c|c|c|c|c|c|c|}
    \hline
   &  & $c_{[0,0]}$ & $c_{[0,1]}$ & $c_{[1,0]}$ & $c_{[1,1]}$ & $c_{[2,0]}$ & $c_{[2,1]}$ \\
\hline
  $id$ &  123 & 1 & $x_2$ & $x_1$ & $x_1x_2$ & $x_1^2$ & $x_1^2x_2$ \\
\hline
  $s_1$ &  213 & 1 & $x_1$ & $x_2$ & $x_2x_1$ & $x_2^2$ & $x_2^2x_1$ \\
\hline
  $s_2$ &  132 & 1 & $x_3$ & $x_1$ & $x_1x_3$ & $x_1^2$ & $x_1^2x_3$ \\
\hline
  $s_1s_2$ &  231 & 1 & $x_3$ & $x_2$ & $x_2x_3$ & $x_2^2$ & $x_2^2x_3$ \\
\hline
  $s_2s_1$ &  312 & 1 & $x_1$ & $x_3$ & $x_3x_1$ & $x_3^2$ & $x_3^2x_1$ \\
\hline
  $s_1s_2s_1$ &  321 & 1 & $x_2$ & $x_3$ & $x_3x_2$ & $x_3^2$ & $x_3^2x_2$ \\
\hline
  \end{tabular}
\\[3pt]
\caption{Invariant classes on $W(A_2)$}
\label{tbl:inv_class_S3}
\end{table}
Repeating the procedure further, we get the following result.
\begin{theorem} Let
$$\mathcal{A}_n= \{ I = [i_1,\ldots, i_{n}]  \; | \;
0 \leqslant i_1 \leqslant n, 0 \leqslant i_2 \leqslant n-1, \ldots , 0 \leqslant i_{n} \leqslant 1\}\; .$$
Then
$$\{ c_I = c_T(\textbf{x}^I) \; | \; I \in \mathcal{A}_n\}$$
is an  $\SS(\ft^*)-$module basis of $H_{\alpha}^*(A_n)$, consisting of invariant classes.
\end{theorem}

By Theorem~\ref{prop:bases_over_invariants} it follows that,
in type $A_n$, $\{ \textbf{x}^I \; | \; I \in \mathcal{A}_n\}$
is a basis of $\SS(\ft^*)$ as an $\SS(\ft^*)^W-$module. Observe that
the top degree class is generated by the top degree Schubert polynomial.

\subsection{Type B}
The set of simple roots of $B_n$ (for $n \geqslant 2$) is $\Delta_0 = \{\alpha_1,\ldots,\alpha_n\}$,
where $\alpha_i = x_i - x_{i+1}$, for $i=1,\ldots,n-1$ and $\alpha_n = x_n$. The set of positive roots is
$$\Delta^{+} = \{ x_i \; | \; 1 \leqslant i \leqslant n \} \cup
\{ x_i \pm  x_j \; | \; 1 \leqslant i < j \leqslant n\} \; .$$
If $S=\{ \alpha_2,\ldots, \alpha_n\}$, then
$$\langle S \rangle = \{ x_i \; | \; 2 \leqslant i \leqslant n \} \cup
\{ x_i \pm  x_j \; | \; 2 \leqslant i < j \leqslant n\}$$
is the set of positive roots for a root system of type $B_{n-1}$, and
$$\Delta^+ \setminus \langle S \rangle = \{ \beta_1=x_1 \} \cup
\{ \beta_j^\pm=x_1\mp x_j \; | \; 2\leqslant j \leqslant n\} \; .$$

Let
\begin{align*}
  \omega_1^+ = & [id] \quad , \quad \omega_1^- = [s_{\beta_1}] \\
  \omega_j^+ = & [s_{\beta_j^+}] = [s_{x_1-x_j}]  \text{ for } 2 \leqslant j \leqslant n \\
  \omega_j^- = &  [s_{\beta_j^-}] = [s_{x_1+x_j}]  \text{ for } 2 \leqslant j \leqslant n\; .
\end{align*}

Then $W/W_S = \{ \omega_1^+, \omega_1^-, \ldots, \omega_n^+, \omega_n^-\}$, and the
graph structure of $W/W_S$ is that of a complete graph with $2n$ vertices. If $\tau$
is the map $\tau \colon W/W_S \to \ft^*$, $\tau(\omega_j^{\epsilon}) = \epsilon x_j$,
with $1 \leqslant j \leqslant n$ and $\epsilon \in \{ +, -\}$, then the axial function $\alpha$
\begin{align*}
  \alpha(\omega_i^{\epsilon_i},\omega_j^{\epsilon_j}) = &
\tau(\omega_i^{\epsilon_i}) - \tau(\omega_j^{\epsilon_j}),  \text{ for } 1 \leqslant i \neq j \leqslant n \\
  \alpha(\omega_i^{\epsilon_i},\omega_i^{-\epsilon_i}) = &
\frac{1}{2} (\tau(\omega_i^{\epsilon_i}) - \tau(\omega_i^{-\epsilon_i})) &
\text{  for } 1 \leqslant i \leqslant n\; .
\end{align*}
Note that although $W/W_S$ and $K_{2n}$ are isomorphic as graphs,
they are not isomorphic as GKM graphs.
One way to see that is to notice that
$$\alpha(\omega_1^+, \omega_1^-) + \alpha(\omega_1^-, \omega_2^-) + \alpha(\omega_2^-, \omega_1^+) =
-x_1 \neq 0 \; .$$
Nevertheless, as in the $K_{2n}$ case, the set of classes $\{1, \tau, \ldots, \tau^{2n-1}\}$
is a basis
for the free $\SS(\ft^*)-$module $H_{\alpha}^*(W/W_S)$.

\begin{figure}[h]
  \psfrag{x1-x2}{\small{$x_1\!-\!x_2$}}
  \psfrag{x1+x2}{\small{$x_1\!+\!x_2$}}
  \psfrag{x1}{\small{$x_1$}}
  \psfrag{x2}{\small{$x_2$}}
  \psfrag{omega1+}{$\omega_1^+$}
  \psfrag{omega1-}{$\omega_1^-$}
  \psfrag{omega2+}{$\omega_2^+$}
  \psfrag{omega2-}{$\omega_2^-$}
  \psfrag{pi}{$\pi$}
  \psfrag{pi(omega1+)}{$\pi^{-1}(\omega_1^+)$}
  \psfrag{pi(omega2+)}{$\pi^{-1}(\omega_2^+)$}
  \psfrag{pi(omega1-)}{$\pi^{-1}(\omega_1^-)$}
  \psfrag{pi(omega2-)}{$\pi^{-1}(\omega_2^-)$}
  \psfrag{id}{id}
  \psfrag{s1}{$s_1$}
  \psfrag{s2}{$s_2$}
  \psfrag{s1s2}{$s_1s_2$}
  \psfrag{s2s1}{$s_2s_1$}
  \psfrag{s1s2s1}{$s_1s_2s_1$}
  \psfrag{s2s1s2}{$s_2s_1s_2$}
  \psfrag{s1s2s1s2}{$s_1s_2s_1s_2=s_2s_1s_2s_1$}
  \includegraphics[height=2in]{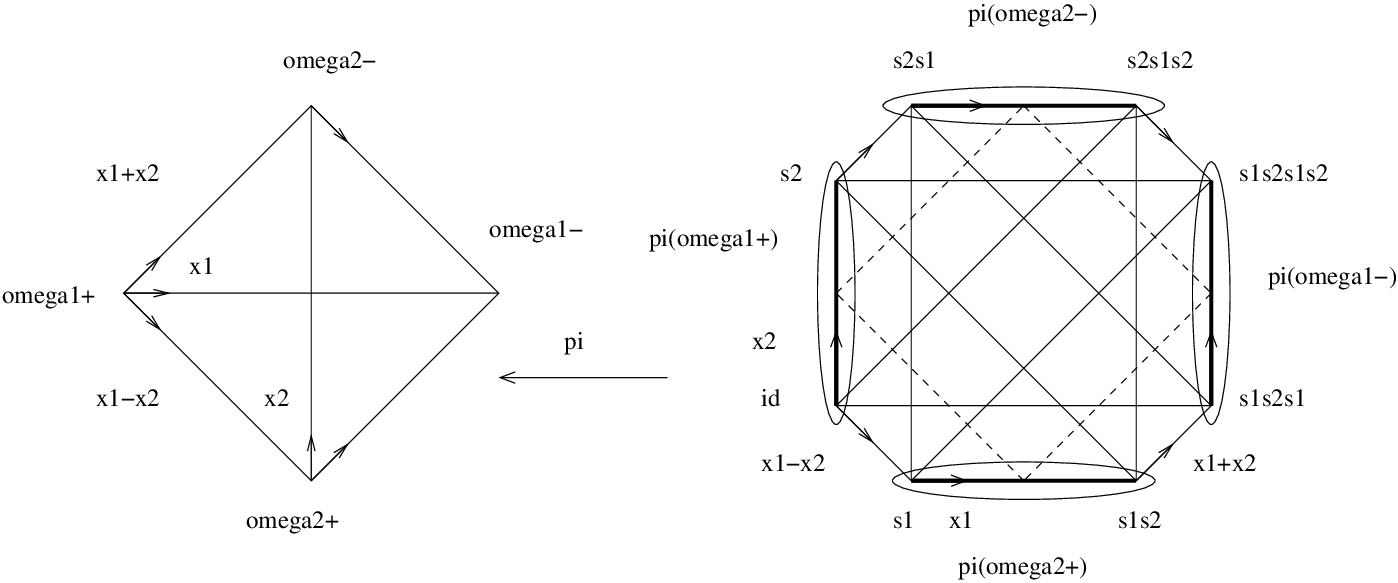}
  \caption{Fibration of $B_2$}
\end{figure}

An alternative description of the Weyl group $W$ is that of the group of signed
permutations $(u,\epsilon)$, with $u \in S_n$ and
$\epsilon = (\epsilon_1,\ldots,\epsilon_n)$, $\epsilon_j = \pm 1$.
The element $(u,\epsilon)$ is represented as $(\epsilon_1 u(1),\ldots,\epsilon_nu(n))$
or by underlying the negative entries.

Then $s_{x_i}$ is just a change of the sign of $x_i$, $s_{x_i-x_j}$
corresponds to the transposition $(i,j)$, with no sign changes, and $s_{x_i+x_j}$
corresponds to the transposition $(i,j)$ with both signs changed. In particular,
$id$ is the identity permutation with no sign changes, $s_{\beta_1}$ is the identity
permutation with the sign of 1 changed, $s_{\beta_j^+}$ is the transposition $(1,j)$
with no sign changes, and $s_{\beta_j^-}$ is the transposition $(1,j)$ with sign
changes for 1 and $j$. In general, if $u \in S_n$ and
$\epsilon = (\epsilon_1,\ldots,\epsilon_n) \in \ZZ_2^n$, then the
element $w=(u, \epsilon) \in W$ acts by $(u, \epsilon) \cdot x_k = \epsilon_k x_{u(k)}$.
Then $W/W_S$ can be identified with $\{ \pm 1, \pm 2, \ldots , \pm n\}$
by $\omega_j^{\epsilon} \to \epsilon j$, and the projection $\pi \colon W \to W/W_S$
is $\pi((u, \epsilon)) = \epsilon_1u(1)$.

For $I=[i_1,\ldots,i_{n}]$, let $c_I =c_T(\textbf{x}^I)  \colon W \to \SS(\ft^*)$ be given by
$$c_I ((u,\epsilon)) = (\epsilon_1x_{u(1)})^{i_1} \dotsb (\epsilon_{n}x_{u(n)})^{i_{n}}\; .$$
Then $c_I \in (H_{\alpha}^*(W))^W$ is an invariant class, and we will construct a basis of the
 free $\SS(\ft^*)-$module $H_{\alpha}^*(W)$ consisting of classes of the type $c_I$, for specific
indices $I$.

When $n=2$, the fiber over 2 is $\pi^{-1}(2) = \{ (2,1), (2,-1)\}$ and is identified with
$W_S=S_2 = \{ 1,-1\}$. A basis for $H_{\alpha}^*(W_S)$ is given by the invariant classes
$\{ c_{[0]}, c_{[1]} \}$, where $c_{[0]} \equiv 1$ and $c_{[1]}(1) = x_1$, $c_{[1]}(-1) = -x_1$.
These classes are extended to the invariant classes $c_{[0,0]}$ and $c_{[0,1]}$ on $W$.

The classes $1$, $\tau$, $\tau^2$, and $\tau^3$ on the base lift to the basic classes
$c_{[0,0]}$, $c_{[1,0]}$, $c_{[2,0]}$, and $c_{[3,0]}$ on $W$. Then a basis for the
free $\SS(\ft^*)-$module $H_{\alpha}^*(W)$ is
$$\{ c_I \; | \; I=[i_1,i_2], \; 0 \leqslant i_1 \leqslant 3, 0 \leqslant i_2 \leqslant 1\} \; .$$

The values of these classes on the elements of $W(B_2)$ are shown in Table~\ref{tbl:B2}.
\begin{table}[ht]
\begin{tabular}{|c|c|c|c|c|c|c|c|c|c|}
\hline
 & & $c_{[0,0]}$ & $c_{[0,1]}$ & $c_{[1,0]}$ & $c_{[1,1]}$ & $c_{[2,0]}$ & $c_{[2,1]}$ & $c_{[3,0]}$ & $c_{[3,1]}$ \\
 \hline
 $id$ & $12$ & $1$ & $x_2$ & $x_1$ & $x_1x_2$ & $x_1^2$ & $x_1^2x_2$ & $x_1^3$ & $x_1^3x_2$ \\
\hline
$s_1$ & $21$ & $1$ & $x_1$ & $x_2$ & $x_1x_2$ & $x_2^2$ & $x_2^2x_1$ & $x_2^3$ & $x_2^3x_1$ \\
\hline
$s_2$ & $1\underline{2}$ & $1$ & $-x_2$ & $x_1$ & $-x_1x_2$ & $x_1^2$ & $-x_1^2x_2$ & $x_1^3$ & $-x_1^3x_2$ \\
\hline
$s_1s_2$ & $2\underline{1}$ & $1$ & $-x_1$ & $x_2$ & $-x_1x_2$ & $x_2^2$ & $-x_2^2x_1$ & $x_2^3$ & $-x_2^3x_1$ \\
\hline
$s_2s_1$ & $\underline{2}1$ & $1$ & $x_1$ & $-x_2$ & $-x_1x_2$ & $x_2^2$ & $x_2^2x_1$ & $-x_2^3$ & $-x_2^3x_1$ \\
\hline
$s_1s_2s_1$ & $\underline{1}2$ & $1$ & $x_2$ & $-x_1$ & $-x_1x_2$ & $-x_1^2$ & $x_1^2x_2$ & $-x_1^3$ & $-x_1^3x_2$ \\
\hline
$s_2s_1s_2$ & $\underline{2}\underline{1}$ & $1$ & $-x_1$ & $-x_2$ & $x_1x_2$ & $x_2^2$ & $-x_2^2x_1$ & $-x_2^3$ & $x_2^3x_1$ \\
\hline
$s_1s_2s_1s_2$ & $\underline{1}\underline{2}$ & $1$ & $-x_2$ & $-x_1$ & $x_1x_2$ & $x_1^2$ & $-x_1^2x_2$ & $-x_1^3$ & $x_1^3x_2$ \\
\hline
\end{tabular}
\\[3pt]
\caption{Invariant classes on $W(B_2)$ }
\label{tbl:B2}
\end{table}

Repeating the procedure further, we get the following result.

\begin{theorem} Let
$$\mathcal{B}_n = \{ I = [i_1,\ldots,i_n] \; | \;
0 \leqslant i_1 \leqslant 2n-1, 0 \leqslant i_2 \leqslant 2n-3, \ldots ,
0 \leqslant i_{n} \leqslant 1\}$$
Then
$$\{ c_I \; |\; I \in \mathcal{B}_n\}$$
is an $\SS(\ft^*)-$module basis of $H_{\alpha}^*(W(B_n))$ consisting of $W-$invariant classes.
\end{theorem}

By Theorem~\ref{prop:bases_over_invariants} it follows that, in type $B_n$,
$\{ \textbf{x}^I \; | \; I \in \mathcal{B}_n\}$ is a basis of $\SS(\ft^*)$
as an $\SS(\ft^*)^W-$module.

\subsection{Type C}
The set of simple roots of $C_n$ (for $n \geqslant 2$) is $\Delta_0 = \{\alpha_1,\ldots,\alpha_n\}$,
where $\alpha_i = x_i - x_{i+1}$, for $i=1,\ldots,n-1$ and $\alpha_n = 2x_n$. The set of positive roots is
$$\Delta^{+} = \{ 2x_i \; | \; 1 \leqslant i \leqslant n \} \cup \{ x_i \pm  x_j \; | \; 1 \leqslant i < j \leqslant n\} \; .$$
If $S=\{ \alpha_2,\ldots, \alpha_n\}$, then
$$\langle S \rangle = \{ 2x_i \; | \; 2 \leqslant i \leqslant n \} \cup \{ x_i \pm  x_j \; | \; 2 \leqslant i < j \leqslant n\}$$
is the set of positive roots for a root system of type $C_{n-1}$, and
$$\Delta^+ \setminus \langle S \rangle = \{ \beta_1=2x_1 \} \cup \{ \beta_j^\pm=x_1\mp x_j \; | \; 2\leqslant j \leqslant n\} \; .$$

Let
\begin{align*}
  \omega_1^+ = & [id] \quad , \quad \omega_1^- =  [s_{\beta_1}] \\
  \omega_j^+ = & [s_{\beta_j^+}] = [s_{x_1-x_j}] \text{ for } 2 \leqslant j \leqslant n \\
  \omega_j^- = & [s_{\beta_j^-}] = [s_{x_1+x_j}] \text{ for } 2 \leqslant j \leqslant n\; .
\end{align*}
This is essentially the same as the type $B$ case, and $W(C_n) \simeq W(B_n)$ is the group of signed
permutations of $n$ letters. Then $W/W_S = \{ \omega_1^+, \omega_1^-, \ldots, \omega_n^+, \omega_n^-\}$,
and the graph structure of $W/W_S$ is that of a complete graph with $2n$ vertices. The axial function
on $W/W_S$ is given by
$$\alpha(\omega_i^{\epsilon_i}, \omega_j^{\epsilon_j}) = \tau(\omega_i^{\epsilon_i}) - \tau(\omega_j^{\epsilon_j})\; ,$$
hence $W/W_S$ is isomorphic, as a GKM graph, with a projection of the complete graph $K_{2n}$.
Then $H_{\alpha}^*(W(C_n)) \simeq H_{\alpha}^*(W(B_n))$, with $\mathcal{B}(C_n) = \mathcal{B}(B_n)$
as a basis consisting of invariant classes.

\subsection{Type D} The set of simple roots of $D_n$ (for $n \geqslant 3$) is
$\Delta_0 = \{\alpha_1,\ldots,\alpha_n\}$, where $\alpha_i = x_i - x_{i+1}$,
for $i=1,\ldots,n-1$ and $\alpha_n = x_{n-1} + x_n$. The set of positive roots is
$$\Delta^{+} = \{ x_i-x_j \; | \; 1 \leqslant i<j \leqslant n \} \cup \{ x_i +  x_j \; | \; 1 \leqslant i < j \leqslant n\} \; .$$
If $S=\{ \alpha_2,\ldots, \alpha_n\}$, then
$$\langle S \rangle = \{ x_i-x_j \; | \; 2 \leqslant i<j \leqslant n \} \cup \{ x_i +  x_j \; | \; 2 \leqslant i < j \leqslant n\} \; . $$
If $n \geqslant 4$, then $\langle S \rangle$ is the set of positive roots for a root system
of type $D_{n-1}$ and if $n=3$, then $\langle S \rangle$ is the set of positive roots of
$A_1 \times A_1$. In both cases
$$\Delta^+ \setminus \langle S \rangle = \{ \beta_i^+ = x_1-x_i \; | \; 2 \leqslant i \leqslant n \} \cup \{ \beta_i^- = x_1 +  x_i \; | \; 2 \leqslant i \leqslant n\} \; .$$

Let
\begin{align*}
  \omega_1^+ & = [id], & \omega_1^- & = [s_{\beta_j^-}s_{\beta_j^+}] = [s_{\beta_j^+}s_{\beta_j^-}], & \text{ for all } 2 \leqslant j \leqslant n \\
  \omega_i^+ & = [s_{\beta_i^+}],  &  \omega_i^- & = [s_{\beta_i^-}], & \text{ for all } 2 \leqslant i \leqslant n \; .
\end{align*}

Then $W/W_S = \{\omega_1^+, \omega_1^-, \ldots, \omega_n^+, \omega_n^- \}$ and the graph structure
of $W/W_S$ is that of the complete $n-$partite graph $K_2^n$, with
partition classes $\{\omega_i^+, \omega_i^-\}$ for $1 \leqslant i \leqslant n$.
If $\tau \colon W/W_S \to \ft^*$ is given by $\tau(\omega_i^{\epsilon}) = \epsilon x_i$,
where $\epsilon \in \{ +, -\}$, then the axial function $\alpha$ on $W/W_S$ is
$$\alpha(\omega_i^{\epsilon_i}, \omega_j^{\epsilon_j}) =
\tau(\omega_i^{\epsilon_i}) - \tau(\omega_j^{\epsilon_j}) =
\epsilon_i x_i - \epsilon_j x_j \; .$$

Then $H_{\alpha}^*(W/W_S)$ is a free $\SS(\ft^*)-$module, and a
Vandermonde determinant argument shows that a basis is given by
$1$, $\tau$, \ldots, $\tau^{2n-2}$, and $\eta = x_1\dotsb x_n \tau^{-1}$.

An alternative description of the Weyl group $W$ is that of the group of
signed permutations $(u,\epsilon)$ with an even number of sign changes.
Then $s_{x_i-x_j}$ corresponds to the transposition $(i,j)$, with no sign changes,
and $s_{x_i+x_j}$ corresponds to the transposition $(i,j)$ with both signs changed.
In particular, $id$ is the identity permutation with no sign changes, $s_{\beta_j^+}$ is
the transposition $(1,j)$ with no sign changes, $s_{\beta_j^-}$ is the transposition $(1,j)$
with sign changes for 1 and $j$, and $s_{\beta_j^+}s_{\beta_j^-}$ is the identity permutation
with the sign changes for $1$ and $j$. In general, if $u \in S_n$ and
$\epsilon = (\epsilon_1,\ldots,\epsilon_n) \in \ZZ_2^n$ with $\epsilon_1 \dotsb \epsilon_n = 1$,
then the element $w=(u, \epsilon) \in W$ acts by $(u, \epsilon) \cdot x_k = \epsilon_k x_{u(k)}$.
Then $W/W_S$ can be identified with $\{ \pm 1, \pm 2, \ldots , \pm n\}$ by
$\omega_i^{\epsilon} \to \epsilon i$, and the projection $\pi \colon W \to W/W_S$
is $\pi((u, \epsilon)) = \epsilon_1u(1)$.

For $I=[i_1,\ldots,i_{n}]$, let $c_I=c_T(\textbf{x}^I) \colon W \to \SS(\ft^*)$ be given by
$$c_I ((u,\epsilon)) = (\epsilon_1x_{u(1)})^{i_1} \dotsb (\epsilon_{n}x_{u(n)})^{i_{n}}\; .$$
Then $c_I \in (H_{\alpha}^*(W))^W$ is an invariant class, and we will construct a
basis of the free $\SS(\ft^*)-$module $H_{\alpha}^*(W)$ consisting of classes of the type $c_I$,
for specific indices $I$.

When $n=3$, the fiber $\pi^{-1}(3)$ of the GKM fiber bundle $\pi \colon D_3 \to K_2^3$ is
$$\pi^{-1}(3) = \{ (3,1,2), (3,2,1), (3,-2,-1), (3,-1,-2)\}$$
and is identified with $W_S=  S_2 \times S_2 =\{ (1,2), (2,1), (-2,-1), (-1,-2)\}$.
Then  $H_{\alpha}^*(W_S)$ is generated by the $W_S-$invariant classes
$\{ c_I \; | \;I \in  \mathcal{D}_2 \}$, where
$$\mathcal{D}_2 = \{ [0,0],[1,0],[2,0],[0,1] \} \; .$$
The classes $1$, $\tau$, $\tau^2$, $\tau^3$, $\tau^4$, $\eta$ on $K_2^3$
lift to the basic classes $c_{[0,0,0]}$, $c_{[1,0,0]}$, $c_{[2,0,0]}$, $c_{[3,0,0]}$,
$c_{[4,0,0]}$, and $c_{[0,1,1]}$.
Then a basis for the free $\SS(\ft^*)-$module $H_{\alpha}^*(W)$ is
$$\{ c_I \; | \; I=[i_1,i_2,i_3] \in \mathcal{D}_3 \}\; ,$$
where $\mathcal{D}_3$ is the set of triples $[i_1,i_2,i_3] \in \ZZ_{\geqslant 0}^3$,
such that $i_1i_2i_3=0$ and either $i_1 \leqslant 4, i_2 \leqslant 2, i_3 \leqslant 1$
or $[i_1,i_2,i_3] =[0,1,2]$ or $[0,3,1]$.

Repeating this process further, we get the following general result.

\begin{theorem} Let $\mathcal{D}_n$ be a set of multi-indices defined inductively by
\begin{enumerate}
 \item $\mathcal{D}_2 = \{[0,0],[1,0],[2,0],[0,1] \}$;
\item $[i_1,\ldots,i_n] \in \mathcal{D}_n$ if
\begin{itemize}
 \item $0 \leqslant i_1 \leqslant 2n-2$ and $[i_2,\ldots,i_n] \in \mathcal{D}_{n-1}$, or
  \item $i_1=0$ and $[i_2-1,\ldots,i_n-1] \in \mathcal{D}_{n-1}$.
\end{itemize}
\end{enumerate}

Then
$$\{ c_I \; | \; I \in \mathcal{D}_n \} \; .$$
is an $\SS(\ft^*)-$module basis of $H_{\alpha}^*(D_n)$ consisting of $W-$invariant classes.
\end{theorem}

By Theorem~\ref{prop:bases_over_invariants} it follows that, in type $D_n$,
$\{ \textbf{x}^I \; | \; I \in \mathcal{D}_n\}$ is a basis of $\SS(\ft^*)$
as a free $\SS(\ft^*)^W-$module.

\section{Symmetrization of Schubert Classes}
\label{sec:symmetrization}
In Section~\ref{sec:classical} we constructed invariant classes for classical groups by
iterating the GKM fiber bundle construction. In this section we present a different method
of constructing invariant classes.

\subsection{Symmetrization of Classes}
\label{sec:sym_classes}

Recall that the ring $H_{\alpha}^*(W)$ consists of the maps
$f \colon W \to \SS(\ft^*)$ such that
$$f(ws_{\beta}) - f(w) \in (w\beta) \SS(\ft^*)$$
for every $w \in W$ and $\beta \in \Delta^+$, and the holonomy action of the Weyl group $W$ is
$$w \cdot f = f^{w} \colon W \to \SS(\ft^*)\; , \quad f^{w}(v) = w^{-1} f(wv) \; .$$

For every $u \in W$, there exists a unique class $\tau_u \in H_{\alpha}^*(W)$, called the
equivariant Schubert class of $u$, that satisfies the following conditions:
\begin{enumerate}
  \item $\tau_u$ is homogeneous of degree $2\ell(u)$, where $\ell(u)$ is the length of $u$;
  \item $\tau_u$ is supported on $\{ v \; | u \preccurlyeq v\}$, where $\preccurlyeq$ is the
  strong Bruhat order, and
  \item $\tau_u$ is normalized by the condition
$$\tau_u(u) = \prod \{\beta \; | \; \beta \in \Delta^+ , u^{-1}\beta \in \Delta^- \}$$
\end{enumerate}

The set $\{\tau_u \; | \; u \in W \}$ of equivariant Schubert classes is a basis of the
$\SS(\ft^*)-$module $H_{\alpha}^*(W)$; however, these classes are not invariant under
the action of $W$ on $H_{\alpha}^*(W)$.

For $f \in H_{\alpha}^*(W)$ we define the $W-$invariant class $f^{sym} \colon W \to \SS(\ft^*)$ by
$$f^{sym} = \frac{1}{|W|} \sum_{w \in W} f^w \; ,$$
where the permuted class $f^{w} \colon W \to \SS(\ft^*)$ is given by $f^w(u) = w^{-1} \cdot f(wu)$, $u \in W$.

For every $w \in W$, the permuted classes $\{\tau_u^w \; | \; u \in W \}$ form a basis of
the $\SS(\ft^*)-$module $H_{\alpha}^*(W)$. The main result of this section is that the
\emph{symmetrized} classes also form a basis of the $H_{\alpha}^*(W)$, and these classes
\emph{are} $W-$invariant.

\subsection{NilCoxeter Rings}
\label{ssec:niHecke}
We start by recalling a few things about nilCoxeter rings. More details are available,
for example, in \cite{Ku}.

These rings are defined for
general Coxeter groups, but we will only need them for Weyl groups, for which
we will use the notation introduced in Section~\ref{sec:flags_as_GKM}.

Let $W$ be a Weyl group, with simple positive roots $\{\alpha_1, \ldots, \alpha_n\}$ and
let $s_i = s_{\alpha_i}$ be the reflection generated by the simple root $\alpha_i$, for
$1 \leqslant i \leqslant n$. The nilCoxeter ring $\mathcal{H}$ is the ring with
generators $\{u_i \; | \; i=1,\ldots,n\}$ satisfying $u_i^2 = 0$ for all
$i=1,\ldots,n$ and the same commutation relations as $\{ s_i \; |\; i=1,\ldots,n\}$.

If $w=s_{i_1}\dotsb s_{i_r}$ is a reduced decomposition of $w \in W$
(hence $\ell(w) = r$), we define
$$u_w = u_{i_1} \dotsb u_{i_r} \; .$$
The definition does not depend on the reduced decomposition, and
$$u_w u_v = \left\{
\begin{array}{ll}
 u_{wv}, & \text{ if } \ell(wv) = \ell(w) + \ell(v) \\
 0, & \text{ otherwise.}
\end{array}
 \right.$$

For every $i=1,\ldots,n$, let $h_i(x) = 1 + xu_i$, where $x$ is a variable
that commutes with all
the generators $u_1, \ldots, u_n$. Then $h_i(x)$ is invertible and $h_i(x)^{-1} = h_i(-x)$.

If $w=s_{i_1}\dotsb s_{i_r}$ is a reduced decomposition of $w \in W$, define
$H_w \in \mathcal{H}\otimes \SS(\ft^*)$ by
\begin{equation}\label{eq:Hw-expl}
\begin{aligned}
H_w =&  h_{i_1}(\alpha_{i_1}) h_{i_2}(s_{i_1}\alpha_{i_2}) \dotsb
h_{i_r} ( s_{i_1} \dotsb s_{i_{r-1}}\alpha_{i_r}) = \\
= & (1+\alpha_{i_1}u_{i_1})(1+s_{i_1}\alpha_{i_2} u_{i_2})\dotsb (1+s_{i_1}\dotsb s_{i_{r-1}}\alpha_{i_r} u_{i_r})
\end{aligned}
\end{equation}
The definition of $H_w$ does not depend on the reduced decomposition of $w$.

In \cite[Theorem 3]{Bi}, Billey showed that
\begin{equation}\label{eq:positive_formula}
  H_w = \sum_{v \in W} \tau_{v}(w) u_v \;
\end{equation}
and used this formula to prove an explicit positive
formula for $\tau_v(w)$, as a sum of products of
positive roots (see also \cite[Appendix D]{AJS}). In particular,
$$\tau_v (w) \in \ZZ_{\geqslant 0}^{\ell(v)}[\alpha_1, \ldots, \alpha_n]$$
is a homogeneous polynomial of degree $\ell(v)$ in the simple positive
roots $\alpha_1,~\ldots,~\alpha_n$, with nonnegative integer coefficients.
Moreover, $H_w$ is invertible, and
\begin{equation}\label{eq:alternating_formula}
  H_w^{-1} = h_{i_r} ( - s_{i_1} \dotsb s_{i_{r-1}}\alpha_{i_r}) \dotsb h_{i_1}( - \alpha_{i_1}) = \sum_{v \in W} (-1)^{\ell(v^{-1})} \tau_{v^{-1}}(w) u_{v} \; .
\end{equation}

\begin{lemma}
\label{lem:main_nilHecke}
  If $w,v \in W$, then
\begin{equation}\label{eq:key_formula}
  H_{wv} = H_w \cdot wH_v \; .
\end{equation}
\end{lemma}

\begin{proof}
If $\ell(v) = 0$, then $v=1$, $H_v = 1$, and the formula is clearly true.

The proof is made in four steps.

\emph{Step 1:}  $v=s_i$ and $\ell(ws_i) = \ell(w) + 1$. Let $w=s_{i_1} \dotsb s_{i_r}$
be a reduced decomposition of $w$; then $ws_i = s_{i_1} \dotsb s_{i_r}s_i$ is a reduced decomposition for $ws_i$, hence
\begin{align*}
  H_{ws_i} = & h_{i_1}(\alpha_{i_1}) h_{i_2}(s_{i_1}\alpha_{i_2}) \dotsb h_{i_r} ( s_{i_1} \dotsb s_{i_{r-1}}\alpha_{i_r})h_{i_r} ( s_{i_1} \dotsb s_{i_{r}}\alpha_{i}) = \\
 = & H_w \cdot h_i( w\alpha_i) = H_w \cdot w h_i(\alpha_i) = H_w \cdot wH_{s_i} \; .
\end{align*}

\emph{Step 2:} $\ell(wv) = \ell(w) + \ell(v)$. If $v = s_{i_1}\dotsb s_{i_r}$ is a reduced
decomposition for $v$, then $ws_{i_1} \dotsb s_{i_k}$ is a reduced decomposition for every
$k=1,\ldots,r$, and hence Step~1 applies in all those cases. Hence
\begin{align*}
  H_{wv} = & H_{ws_{i_1}\dotsb s_{i_{r-1}} s_{i_r}} = H_{ws_{i_1}\dotsb s_{i_{r-1}}} \cdot ws_{i_1}\dotsb s_{i_{r-1}} h_{i_r}(\alpha_{i_r}) = \\
= & H_{ws_{i_1}\dotsb s_{i_{r-1}}} \cdot w h_{s_{i_r}}(s_{i_1}\dotsb s_{i_{r-1}}\alpha_{i_r}) = \\
 = & H_w \cdot wh_{i_1}(\alpha_{i_1}) \dotsb w h_{i_r}(s_{i_1}\dotsb s_{i_{r-1}}\alpha_{i_r}) =
 H_w \cdot wH_v \; .
\end{align*}

\emph{Step 3:} $v=s_i$ and $\ell(ws_i) = \ell(w) - 1$. Let $w=s_{i_1} \dotsb s_{i_r}$ be a
reduced decomposition of $w$; then by the Exchange Condition, there exists an index $k$ such
that $ws_i = w_1 w_2$, where $w_1 = s_{i_1} \dotsb s_{i_{k-1}}$ and
$w_2 = s_{i_{k+1}} \dotsb s_{i_r}$ are reduced decompositions. Let $j=i_k$. Then
$w=w_1 s_j w_2$, $s_j w_2 = w_2 s_i$, and $\ell(w_2s_i) = \ell(w_2) +1$.

Then, by the result of Step 2, we have $H_w = H_{w_1 s_jw_2} = H_{w_1} \cdot w_1 H_{s_jw_2}$,
hence
\begin{align*}
  H_w \cdot wH_{s_i} = & H_{w_1} \cdot w_1 H_{s_jw_2} \cdot w_1s_jw_wH_{s_i} = H_{w_1} \cdot w_1 (H_{w_2s_i} \cdot w_2s_iH_{s_i}) = \\
 = & H_{w_1} \cdot w_1 (H_{w_2} \cdot w_2H_{s_i} \cdot w_2s_iH_{s_i}) = H_{w_1} \cdot w_1 (H_{w_2} \cdot w_2 (H_{s_i} \cdot s_iH_{s_i})) \; .
\end{align*}

But $H_{s_i} \cdot s_iH_{s_i} = (1+\alpha_i u_i)(1-\alpha_i u_i) = 1$, hence
$$H_w \cdot wH_{s_i} = H_{w_1} \cdot w_1 H_{w_2} = H_{w_1w_2} = H_{ws_i} \; .$$

At this point we have proved that the formula is true for all $w$ and $v=s_i$.

\emph{Step 4:} For the general case we follow the same argument as for Step 2, using Step 1 or 3
to move over a simple reflection in the reduced decomposition of $v$.
\end{proof}

We use Lemma~\ref{lem:main_nilHecke} to obtain the transition matrices between a basis of
permuted Schubert classes and the original basis of Schubert classes.

\begin{theorem}\label{th:permuted_classes}
Let $a,b,w \in W$. Then
\begin{align}
  \label{eq:Schub_decomposition}
  \tau_a = & \sum_{b \leqslant_{\,L} a} \tau_{ab^{-1}}(w^{-1}) \tau_b^w  \; , \\
\label{eq:permuted_decomposition}
  \tau_a^w = & \sum_{b \leqslant_{\,L} a} (-1)^{\ell(ba^{-1})} \tau_{ba^{-1}}(w^{-1}) \tau_b  \; .
\end{align}
where $\leqslant_{\,L}$ is the left weak order, defined by
$v \leqslant_{\,L} u \Longleftrightarrow \ell(uv^{-1}) = \ell(u) -\ell(v)$.
\end{theorem}

\begin{proof}
  Let $v \in W$. By equation \eqref{eq:key_formula} we have
\begin{equation}
\label{eq:Hv}
  H_{v} = H_{w^{-1}} \cdot w^{-1}H_{wv}
\end{equation}
which, using equation \eqref{eq:positive_formula} and identifying the corresponding coefficients, yields
$$\tau_a(v) = \sum_{\substack{tb=a \\ \ell(t) + \ell(b) = \ell(a)}}
\tau_t(w^{-1}) \cdot w^{-1}\tau_b(wv) = \sum_{b \leqslant_{\,L} a} \tau_{ab^{-1}}(w^{-1})\tau_b^w(v)\, .$$
Since this is true for all $v \in W$, we get \eqref{eq:Schub_decomposition}.

From \eqref{eq:Hv} we get
$$w^{-1}H_{wv} = H_{w^{-1}}^{-1} H_v \; ,$$
which, using \eqref{eq:positive_formula}-\eqref{eq:alternating_formula} and identifying
the corresponding coefficients, yields
$$\tau_a^w(v) = \sum_{\substack{tb=a \\ \ell(t) + \ell(b) = \ell(a)}}
(-1)^{\ell(t^{-1})} \tau_{t^{-1}} (w^{-1}) \tau_b(v) = \sum_{b \leqslant_{\,L} a} (-1)^{\ell(ba^{-1})} \tau_{ba^{-1}}(w^{-1}) \tau_b(v)\; .$$
Since this is true for all $v \in W$, we get \eqref{eq:permuted_decomposition}.
\end{proof}

If $w \in W$ then ${\mathcal B}^w = \{ \tau_u^w \, | \, u \in W\}$ is a basis of
$H_{\alpha}^*(W)$ as an $\SS(\ft^*)-$module. By~\eqref{eq:Schub_decomposition} the
transition matrix $a^w$ between $\mathcal{B}^w$ and the basis ${\mathcal B} = \{ \tau_u \, | \, u \in W\}$
is the lower triangular (with respect to the weak left order) matrix
$$a_{u,v}^w = \left\{ \begin{array}{ll}
   (-1)^{\ell(vu^{-1})} \tau_{vu^{-1}}(w^{-1}) \; , & \; \text{ if } v \leqslant_L u \\
  0 \; , & \; \text{ otherwise.}
\end{array}
\right.$$
Since $\tau_{vu^{-1}}(w^{-1}) \in \ZZ_{\geqslant 0}[\alpha_1,\ldots,\alpha_n]$ is homogeneous of
degree $\ell(vu^{-1})$, we have
$$a_{u,v}^w \in \ZZ_{\geqslant 0}^{\ell(u)-\ell(v)}[-\alpha_1,\ldots,-\alpha_n]\; .$$
Hence the nonzero entries of $a^w$ are homogeneous polynomials in the negative simple roots,
with non-negative integer coefficients, and the diagonal entries are 1.
By~\eqref{eq:permuted_decomposition},
the inverse of $a^w$ is the lower triangular matrix $b^w$ with entries
$$b_{u,v}^w = \left\{ \begin{array}{ll}
   \tau_{uv^{-1}}(w^{-1}) \; , & \; \text{ if } v \leqslant_L u \\
  0 \; , & \; \text{ otherwise.}
\end{array}
\right.$$
The nonzero entries of $b^w$ are homogeneous polynomials in the positive simple roots, with
non-negative integer coefficients, and, again, the diagonal entries are 1.

\subsection{Symmetrized Schubert Classes}
\label{ssec:sym_schubs}

The next result gives the decompositions of symmetrized Schubert classes in terms of Schubert
classes and proves that the set
$\mathcal{B}^{sym} = \{ \tau_u^{sym} \, | \, u \in W\}$ of symmetrized classes is a
basis of $H_{\alpha}^*(W)$.

\begin{theorem}\label{th:symmetrized-classes}
 For every $u \in W$, let $\tau_u^{sym}$ be the symmetrization of $\tau_u$. If
\begin{equation}
  \tau_u^{sym} = \sum_{v \in W} a_{u,v} \tau_v\; ,
\end{equation}
is the decomposition of $\tau_u^{sym}$ in the Schubert basis, then
\begin{enumerate}
  \item The matrix $(a_{u,v})_{u,v}$ is lower triangular with respect to the left weak order:
$$a_{u,v} \neq 0 \Longrightarrow v \leqslant_L u \; .$$
  \item The entries on the diagonal are all 1:
$$a_{u,u} = 1 \quad \text{ for all } u \in W \; .$$
\item The set $\mathcal{B}^{sym} = \{ \tau_u^{sym} \; | \; u \in W\}$
is a basis of the $\SS(\ft^*)-$module $H_{\alpha}^*(W)$.
\end{enumerate}
\end{theorem}

\begin{proof}
If $u \in W$ then
$$\tau_u^{sym} = \frac{1}{|W|} \sum_{w \in W} \tau_u^w = \frac{1}{|W|} \sum_{w \in W} \sum_{v \leqslant_L u} a_{u,v}^w \tau_v = \sum_{v \leqslant_L u} \left( \frac{1}{|W|} \sum_{w \in W} a_{u,v}^w \right) \tau_v \; .$$
Therefore
$$a_{u,v} = \frac{1}{|W|} \sum_{w \in W} a_{u,v}^w \; ,$$
hence $(a_{u,v})_{u,v}$ is lower triangular with respect to the left weak order, with entries
on the diagonal equal to 1. Such a matrix is invertible, and since ${\mathcal B}$ is a basis
of $H_{\alpha}^*(W)$, it follows that ${\mathcal B}^{sym}$ is also a basis.
\end{proof}

\begin{remark}
  For $v \leqslant_{\,L}u$ we have
  $$|W| a_{u,v} \in \ZZ_{\geqslant 0}^{\ell(u)-\ell(v)} [-\alpha_1,\ldots,-\alpha_n] \; ,$$
  because for all $w\in W$, $a_{u,v}^w$ is a homogeneous polynomial  of degree $\ell(u) - \ell(v)$
  in the negative simple roots, with non-negative integer coefficients.
\end{remark}

\subsection{Decomposition of Invariant Classes} Theorem~\ref{th:symmetrized-classes}
gives the decomposition of a symmetrized Schubert class $\tau_{u}^{sym}$ in the
Schubert basis $\{ \tau_w \}_w$. In this section we show how a general invariant
class $c_f = c_T(f) \in H_{\alpha}^*(W)^W$, defined by
\eqref{eq:equiv_charact_hom}, decomposes in the Schubert basis.

For $i=1,\ldots,n$ let $\partial_i \colon \SS(\ft^*) \to \SS(\ft^*)$ be the divided
difference operator
$$\partial_i E = \frac{E-s_i \cdot E}{\alpha_i}\; .$$
If $w=s_{i_1}s_{i_2} \dotsb s_{i_m}$ is a
reduced decomposition for $w\in W$, let $\epsilon(w) = (-1)^{\ell(w)}$ and
$$\partial_w = \partial_{i_1}\partial_{i_2} \dotsb \partial_{i_m}\; ;$$
the notation is justified by the fact that the result of the composition
depends only on $w$ and not on the reduced decomposition of $w$.

\begin{proposition}
\label{prop:inv_decomp}
If $f\in \SS(\ft^*)$, then
  $$c_f = \sum_{w\in W} (\epsilon(w) \partial_w f) \, \tau_w \; .$$
\end{proposition}

\begin{proof}
  We have to show that for every $v \in W$ we have
  $$v\cdot f = \sum_{w \in W} (\epsilon(w) \partial_w f)  \,\tau_w(v)\; ,$$
  and we prove this by induction on the length $\ell(v)$ of $v$.

When $\ell(v) = 0$ we have $v=1$ and the only Schubert class $\tau_w$ that has a
nonzero value at $v=1$ is the one corresponding to $w=1$, with $\tau_1(1) = 1$.
Then $\partial_w f = f$ and the formula is obviously true.

Now suppose the formula is true for all $v$ such that $\ell(v) \leqslant k$ and
let $u\in W$ such that $\ell(u) = k+1$. Then $u$ can be written as $u=s_i v$ for
some $i=1,\ldots,n$ and some $v\in W$ such that $\ell(v) = \ell(u)-1 = k$. Then
$$\sum_{w\in W} (\epsilon(w)\partial_w f)  \tau_w(u) = \sum_{w\in W} (\epsilon(w) \partial_w f) \tau_w(s_iv) \, .$$
But
$$\tau_w(s_iv) = s_i \tau_w(v) + \left\{\begin{array}{ll}\alpha_i s_i\tau_{s_iw}(v), & \text{if } s_iw \prec w \\
0 & \text{ otherwise}
\end{array} \right.$$
This follows from $\tau_w(s_iv) = s_i \cdot \tau_w^{s_i}(v)$ and our formula
for $\tau_w^{s_i}$ or from \cite{Kn}. Hence
$$\sum_{w\in W} (\epsilon(w)\partial_w f)\, \tau_w(s_iv) = \sum_{w\in W} (\epsilon(w)\partial_w f)\, s_i \tau_w(v) + \sum_{s_iw\prec w} (\epsilon(w)\partial_w f)\, \alpha_i s_i\tau_{s_iw}(v)\; .$$
However, since
$$\partial_i \partial_{s_iw} = \left\{ \begin{array}{ll}
  \partial_w & \text{ if } s_iw \prec w \\
  0 & \text{ otherwise} \; ,
\end{array}\right.$$
we can rewrite the last sum and using $\epsilon(w) = -\,\epsilon(s_iw)$ we get
\begin{align*}
  \sum_{w\in W} &(\epsilon(w)\partial_w f) \tau_w(s_iv)  = \sum_{w\in W} (\epsilon(w)\partial_w f) s_i \tau_w(v) - \sum_{w\in W} (\epsilon(s_iw)\partial_i\partial_{s_iw} f) \alpha_i s_i\tau_{s_iw}(v) = \\
  & =  \sum_{w\in W} (\epsilon(w)\partial_w f) s_i \tau_w(v) - \sum_{w\in W} (\epsilon(w)\partial_i\partial_{w} f) \alpha_i s_i\tau_{w}(v) = \\
  & =  \sum_{w\in W} (\epsilon(w)\partial_w f) s_i \tau_w(v) - \sum_{w\in W} \epsilon(w) \frac{\partial_w f - s_i \partial_w f}{\alpha_i} \, \alpha_i s_i\tau_{w}(v) = \\
  & = \sum_{w \in W} \epsilon(w)s_i (\partial_w f) s_i\tau_{w}(v) = s_i \sum_{w\in W} (\epsilon(w)\partial_w f) \tau_w(v)\; .
\end{align*}
From the induction hypothesis the last sum is $v\cdot f$ and therefore
$$\sum_{w\in W} (\epsilon(w)\partial_w f) \tau_w(u) = \sum_{w\in W} (\epsilon(w)\partial_w f) \tau_w(s_iv) = s_i \cdot( v \cdot  f) = (s_i v) \cdot f = u\cdot f\; .$$

The induction is complete and that concludes the proof.
\end{proof}

\begin{remark} Comparing Proposition~\ref{prop:inv_decomp} with \cite[p. 65]{Hi}, we see that
$$c_T \colon \SS(\ft^*) \to H_K^*(M) =  H_{\alpha}^*(W)^W$$
is an equivariant version of the characteristic homomorphism
$c\colon \SS(\ft^*) \to H^*(M)$.
\end{remark}

\subsection{Decomposition of symmetrized Schubert classes}
\label{ssec:symSchubdecomp}

For $w~\in~W$, the symmetrized Schubert class $\tau_{w}^{sym}$
is an invariant class and $\tau_w^{sym} = c_T(f_{w})$, where
$$f_w = \tau_w^{sym}(1) = \frac{1}{|W|} \sum_{v \in W} v^{-1} \cdot \tau_w (v) \in \SS(\ft^*)\; .$$

In this subsection we prove a
simple formula for $f_w$ and we use it to revisit the decomposition
of $\tau_w^{sym}$ in terms of the equivariant classes $\tau_u'$s.

\begin{theorem}\label{thm:fw} Let $w_0$ be the longest element of $W$ and
$\Lambda_0 = \tau_{w_0}(w_0) = \prod_{\alpha \succ 0} \alpha$ the product
of all positive roots. If $w \in W$, then
\begin{equation}\label{eq:fw}
f_w = \frac{\varepsilon(w)}{|W|} \partial_{w^{-1}w_0} (\Lambda_0)
\end{equation}
\end{theorem}

\begin{proof} We prove this result by descending induction on $\ell(w)$.
For $w=w_0$ we have
$$f_{w_0} = \frac{1}{|W|} \sum_{v \in W} v^{-1}\cdot \tau_{w_0}(v) = \frac{1}{|W|} w_0^{-1} \cdot \Lambda_0 = \frac{\varepsilon(w_0)}{|W|} \partial_{w_0^{-1}w_0} (\Lambda_0)\; , $$
because the action of $w_0^{-1} = w_0$ changes all positive roots to negative roots.

Suppose that \eqref{eq:fw} is true for $u$ and let $w = us_i \prec w$ for a
simple reflection $s_i$. Then  $w^{-1}w_0 = s_iu^{-1}w_0$ and $\ell(w) = \ell(u)-1$.
This implies
$$\ell(s_i u^{-1}w_0) = \ell(w_0) - \ell(s_iu{-1}) = \ell(w_0) - \ell(us_i) = \ell(w_0) - \ell(u) + 1 = \ell(s_i)+ \ell(u^{-1}w_0)$$
and therefore $\partial_{w^{-1}w_0} = \partial_{i}\partial_{u^{-1}w_0}\; .$
Hence the right hand side of \eqref{eq:fw} becomes
$$
\frac{\varepsilon(w)}{|W|} \partial_{w^{-1}w_0} (\Lambda_0) =
-\frac{\varepsilon(u)}{|W|} \partial_{i} \partial_{u^{-1}w_0}
(\Lambda_0) = - \partial_i  (f_u) =  - \frac{1}{\alpha_i}
(f_u - s_i \cdot f_u)\; .$$
But
$$ -\frac{1}{\alpha_i} (f_u - s_i \cdot f_u) = \frac{-1}{|W|\alpha_i}
 \left( \sum_{v \in W} v^{-1} \cdot \tau_u(v) - \sum_{v \in W} s_iv^{-1} \cdot \tau_u(v)
\right)$$
and, after a change of variables in the second sum and using \cite[Prop. 2]{Kn},
\begin{align*}
\frac{\varepsilon(w)}{|W|} & \partial_{w^{-1}w_0} (\Lambda_0) =
\frac{-1}{|W|\alpha_i} \sum_{v \in W} v^{-1}\cdot (\tau_u(v) -\tau_u(vs_i)) = \\
= &
\frac{1}{|W|} \sum_{v \in W} v^{-1} \cdot \left( \frac{\tau_u(v) - \tau_u(vs_i)}{v\cdot (-\alpha_i)}\right) =  \frac{1}{|W|} \sum_{v \in W} v^{-1} \cdot \tau_w(v) = f_w \; ,
\end{align*}
completing the proof.
\end{proof}

\begin{remark}
Combining Theorem ~\ref{thm:fw} and Proposition~\ref{prop:inv_decomp} we get
\begin{equation*}
\tau_{w}^{sym} = \frac{1}{|W|} \sum_{v \leqslant_L w} (\varepsilon(v)\varepsilon(w) \partial_{vw^{-1}w_0} (\Lambda_0) )\tau_v \; ,
\end{equation*}
hence the entries of the transition matrix in Theorem~\ref{th:symmetrized-classes}
are given, for $v \leqslant_L u$, by
\begin{equation}
a_{u,v} = \frac{1}{|W|} \varepsilon(u)\varepsilon(v) \partial_{vu^{-1}w_0} (\Lambda_0) = \frac{\varepsilon(vu^{-1})}{|W|} \partial_{vu^{-1}w_0} (\Lambda_0)\; .
\end{equation}
\end{remark}
\begin{remark}
Since $\{ \tau_w^{sym}\}_w$ is a
basis of $H_{\alpha}^*(W)$ over $\SS(\ft^*)$,
Theorem~\ref{prop:bases_over_invariants} implies that
$\{f_w\}_w$ is a basis of $\SS(\ft^*)$
over $\SS(\ft^*)^W$. Therefore, if $\Lambda_0$ is the product
of positive roots, then $\{ \partial_w \Lambda_0\}_w$ is a
basis of $\SS(\ft^*)$ over $\SS(\ft^*)^W$.
\end{remark}

\end{document}